\documentclass[11pt]{amsart}
\usepackage{amssymb,hyperref}
\hypersetup{pdfpagemode=FullScreen,colorlinks=true}
\font\smc=cmcsc9
\font\small=cmr9

\newtheorem{thm}{Theorem}[section]
\newtheorem{sthm}{Theorem}
\newtheorem{scor}[sthm]{Corollary}
\newtheorem{Prop}[thm]{Proposition}
\newtheorem{lemma}[thm]{Lemma}

\newtheorem{co}[thm]{Corollary}
\newtheorem{Def}[thm]{Definition}
\newtheorem{rem}[thm]{Remark}
\newtheorem{ex}[thm]{Example}
\newtheorem{qu}{Question}

\newcommand{\bc}{{\mathbb C}}
\newcommand{\bF}{{\mathbb F}}
\newcommand{\bh}{{\mathbb H}}

\newcommand{\bo}{{\mathbb O}}
\newcommand{\bq}{{\mathbb Q}}
\newcommand{\br}{{\mathbb R}}
\newcommand{\bz}{{\mathbb Z}}
\newcommand{\ra}{\rightarrow}

\newenvironment{pf}{\begin{trivlist}\item[]{\bf Proof:\ }}
{\mbox{}\hfill\rule{.08in}{.08in}\end{trivlist}}
\def\chigg{\chi(\Gamma,G)}

\begin{document}
\title[Density of Zariski density]
{Density of Zariski density for surface groups}
\author{Inkang Kim and Pierre Pansu}
\date{}

\begin{abstract}
We show that a surface group contained in a reductive real algebraic group can be deformed to become Zariski dense, unless its Zariski closure acts transitively on a Hermitian symmetric space of tube type. This is a kind of converse to a rigidity result of Burger, Iozzi and Wienhard.

\medskip\noindent{\smc R\'esum\'e}. {\small On montre qu'un groupe de surface contenu dans un groupe alg\'ebrique r\'eel r\'eductif peut \^etre d\'eform\'e pour devenir Zariski dense, sauf si son adh\'erence de Zariski agit transitivement sur un espace sym\'e\-tri\-que hermitien de type tube. C'est une r\'eciproque partielle d'un th\'e\-o\-r\`eme de rigidit\'e d\^u \`a Burger, Iozzi et Wienhard.}

\end{abstract}

\footnotetext[1]{I. Kim gratefully acknowledges the partial support
of NRF grant  ((R01-2008-000-10052-0) and a warm support of IHES
during his stay.} \footnotetext[2]{P. Pansu, Univ Paris-Sud,
Laboratoire de Math\'ematiques d'Orsay, Orsay, F-91405}
\footnotetext[3]{\hskip 42pt CNRS, Orsay, F-91405.}

\maketitle

\section{Introduction}

We are interested in the dimension of representation varieties. If $G$ is a Lie group and $\Gamma$ a finitely generated group, the corresponding representation variety is the space $Hom(\Gamma,G)$ of homomorphisms $\Gamma\to G$. This space has a {\em virtual dimension} $\mathrm{vdim}_{\phi}(Hom(\Gamma,G))$ at each representation $\phi$, which will be defined in \ref{defvirtual}. For instance, if $G$ is a semisimple real algebraic group and $\Gamma$ is the fundamental group of a closed surface of Euler characteristic $\chi(\Gamma)$, $\mathrm{vdim}_{\phi}(Hom(\Gamma,G))=(1-\chi(\Gamma))\mathrm{dim}(G)$.

In general, we expect $Hom(\Gamma,G)$ to be a stratified space with one open stratum of dimension $\mathrm{vdim}_{\phi}(Hom(\Gamma,G))$ and lower dimensional strata. Examples (see below) show that there may exist open strata of lower dimension. This suggests the following definitions.

\begin{Def}
\label{defrigid}
Let $\Gamma$ be the fundamental group of a closed surface. Let $G$ be a real Lie group. Say a homomorphism $\phi:\Gamma\to G$ is
\begin{itemize}
  \item {\em smooth} if the Zariski tangent space of $Hom(\Gamma,G)$ at $\phi$ is equal to its virtual dimension (this implies that $Hom(\Gamma,G)$ is a smooth manifold near $\phi$ and its dimension equals $\mathrm{vdim}(Hom(\Gamma,G))$;
  \item {\em flexible} if smooth homomorphisms are dense in a neighborhood of $\phi$ in $Hom(\Gamma,G)$;
  \item {\em rigid} if there is a neighborhood of $\phi$ in $Hom(\Gamma,G)$ which contains no smooth homomorphisms.
\end{itemize}
\end{Def}

The sets $\mathcal{S}$ of smooth homomorphisms, $\mathcal{F}$ of flexible homomorphisms and $\mathcal{R}$ of rigid homomorphisms are open in $Hom(\Gamma,G)$. By definition, $\mathcal{F}$ is the interior of the closure of $\mathcal{S}$, and $\mathcal{R}$ is the complement of the closure of $\mathcal{S}$. Roughly speaking, $\mathcal{F}$ is the open stratum of dimension $\mathrm{vdim}(Hom(\Gamma,G))$, and $\mathcal{R}$ is the union of open strata of lower dimensions.

\begin{qu}
\label{question0} When do there exist homomorphisms $\Gamma\to G$
which are rigid? flexible? neither rigid nor flexible?
\end{qu}

This question is related to Zariski density. Indeed (see Corollary \ref{dense}), if $G$ is connected algebraic and if $genus(\Gamma)\geq \mathrm{dim}(G)^2$, the set $\mathcal{Z}$ of homomorphisms $\Gamma\to G$ whose image is Zariski dense is a dense subset of $\mathcal{S}$. Therefore, under the genus restriction, ``smooth'' can be replaced by ``Zariski dense'' in the definitions of
$\mathcal{F}$ and $\mathcal{R}$. In particular, rigid homomorphisms are those which cannot be slightly perturbed to become Zariski dense.

\subsection{A long rigidity story}
\label{rigidity}

We collect examples of rigid homomorphisms.

\subsubsection{Toledo invariants}
Let $X$ be a Hermitian symmetric space, with K\"ah\-ler  form
$\Omega$ (the metric is normalized so that the minimal sectional
curvature equals $-1$). Let $\Sigma$ be a closed surface of negative
Euler characteristic, let $\Gamma=\pi_1 (\Sigma)$ act isometrically
on $X$. Pick a smooth equivariant map $\tilde{f}:\tilde{\Sigma}\to
X$.

\begin{Def}
Define the Toledo invariant of the action $\rho:\Gamma\to Isom(X)$ by
\begin{eqnarray*}
T_{\rho}=\frac{1}{2\pi}\int_{\Sigma}\tilde{f}^{*}\Omega.
\end{eqnarray*}
\end{Def}

Then
\begin{enumerate}
  \item $T_{\rho}$ depends continuously on $\rho$.
  \item There exists $\ell_X \in\bq$ such that $T_{\rho}\in\ell_X \bz$.
  \item $|T_{\rho}|\leq|\chi(\Sigma)|\mathrm{rank}(X)$.
\end{enumerate}

\subsubsection{Rigidity in rank 1}
\begin{ex}
When $X=H^1_{\bc}$ is the unit disk, inequality $|T_{\rho}|\leq|\chi(\Sigma)|$ is due to J. Milnor, \cite{Milnor}. Furthermore $\ell_X =1$, $T$ takes all integer values between $-|\chi(\Sigma)|$ and $|\chi(\Sigma)|$.
\end{ex}

\begin{thm}
\label{PU(1,1)}
{\em (W. Goldman, \cite{G0}).} Let $X=H^1_{\bc}$. The level sets of $T$ coincide with the connected components of the character variety $\chi(\Gamma,PU(1,1))$. Furthermore $|T_{\rho}|=|\chi(\Sigma)|$ if and only if $\rho(\Gamma)$ is discrete and cocompact in $PU(1,1)=Isom(H^1_{\bc})$.
\end{thm}

Note that all components of $\chi(\Gamma,PU(1,1))$ have the same dimension $3|\chi(\Sigma)|$.

\begin{thm}
\label{Toledo}
{\em (D. Toledo, 1979, 1989, \cite{To}).} Let $X=H^n_{\bc}$ have rank 1. Then $|T_{\rho}|\leq|\chi(\Sigma)|$. Furthermore, $|T_{\rho}|=|\chi(\Sigma)|$ if and only if $\rho(\Gamma)$ stabilizes a complex geodesic $H^1_{\bc}$ in $X$ and acts cocompactly on it.
\end{thm}

It follows that, for $n\geq 2$, different components of $\chi(\Gamma,PU(n,1))$ can have different dimensions. Actions with maximal Toledo invariant are rigid. They form an open subset of $\chi(\Gamma,PU(n,1))$ of dimension $-\chi(\Sigma)\mathrm{dim}(P(U(1,1)\times U(n-1)))$. Actions with vanishing Toledo invariant are flexible, they form an open subset of $\chi(\Gamma,PU(n,1))$ of dimension $-\chi(\Sigma)\mathrm{dim}(PU(n,1))$.

\subsubsection{Higher rank}
\begin{Def}
Actions $\rho$ such that $|T_{\rho}|=|\chi(\Sigma)|\mathrm{rank}(X)$ are called {\em maximal representations}.
\end{Def}

\begin{ex}
Pick cocompact actions $\rho_1 ,\ldots,\rho_r$ of $\Gamma$ on $H^1_{\bc}$. Then the direct sum representation on the polydisk $(H^1_{\bc})^r$ is maximal. When the polydisk is embedded in a larger symmetric space of rank $r$, it remains maximal. It follows that all Hermitian symmetric spaces admit maximal representations.
\end{ex}

\begin{Prop}
{\em (Burger, Iozzi, Wienhard, \cite{BIW}).} In case $X$ is Siegel's upper half space (i.e. $Isom(X)=Sp(n,\br)$), some of these actions can be bent to become Zariski dense.
\end{Prop}

But this may fail for other Hermitian symmetric spaces.

\begin{thm}
{\em (L. Hern\`andez Lamoneda, \cite{HL}, S. Bradlow, O. Garc\'{\i}a-Prada, P. Gothen, \cite{BGG}).} Maximal reductive representations of $\Gamma$ to $PU(p,q)$, $p\leq q$, can be conjugated into $P(U(p,p)\times U(q-p))$.
\end{thm}

\subsubsection{Tube type}

The key to these different behaviours of different Hermitian symmetric spaces seems to lie in the notion of tube type domains, as was discovered by M. Burger, A. Iozzi and A. Wienhard.

\begin{Def}
Say a Hermitian symmetric space is of {\em tube type} if it can be realized as a domain in $\bc^n$ of the form $\br^n +iC$ where $C\subset\br^n$ is a proper open cone.
\end{Def}

\begin{ex}
The Grassmannian $\mathcal{D}_{p,q}$, $p\leq q$, with isometry group $PU(p,q)$ is of tube type iff $p=q$.

Products of tube type spaces are of tube type, so polydisks are of tube type.
\end{ex}

\begin{rem}
All maximal tube type subsymmetric spaces in a Hermitian symmetric space are conjugate. For instance, the maximal tube type subsymmetric space in $\mathcal{D}_{p,q}$ is $\mathcal{D}_{p,p}$.
\end{rem}

\begin{thm}
{\em (Burger, Iozzi, Wienhard, \cite{BIW}).} Let $\Gamma$ be a closed surface group and $X$ a Hermitian symmetric space. Every maximal representation $\Gamma\to Isom(X)$ stabilizes a tube type subsymmetric space $Y$. Conversely, for every tube type Hermitian symmetric space $X$, $Isom(X)$ admits Zariski dense maximal surface subgroups.
\end{thm}

In particular, maximal representations of surface groups in non tube type Hermitian symmetric spaces are rigid.

\begin{ex}
In case $X$ is the $n$-ball $\mathcal{D}_{1,n}$ (resp. $\mathcal{D}_{p,q}$), one recovers Toledo's (resp. Bradlow et al.) results.
\end{ex}

\subsection{New flexibility results}

Here are answers to Question \ref{question0}.

First, there always exist flexible homomorphisms.

\begin{sthm}
\label{thm0}
Let $G$ be a connected reductive real algebraic group and let $\Gamma$ be the fundamental group of a closed surface of genus $>1$. Then the trivial homomorphism $\Gamma\to G$ can be deformed into flexible homomorphisms.
\end{sthm}

Second, ``neither flexible nor rigid'' does not show up if genus is large enough.

\begin{sthm}
\label{thm1}
Let $G$ be a connected reductive real algebraic group and let $\Gamma$ be the fundamental group of a closed surface of genus $\geq 2\mathrm{dim}(G)^2$. Then homomorphisms $\Gamma\to G$ are either flexible or rigid.
\end{sthm}
In other words, the sets $\mathcal{F}$ of flexible and $\mathcal{R}$ of rigid homomorphisms are unions of connected components of $Hom(\Gamma,G)$.

\medskip

This leads to a refined question.

\begin{qu}
\label{question1} Which surface groups in real reductive algebraic
groups are flexible?
\end{qu}

As far as the flexibility of a homomorphism $\phi$ is concerned, a key role is played by the center of the centralizer of the image of $\phi$. It splits the Lie algebra of $G$ into real root spaces $\mathfrak{g}_{\lambda,\br}$ which carry natural symplectic structures. In particular, to each pure imaginary root $\lambda$, there corresponds a symplectic representation $\rho_{\lambda}$ and a Toledo invariant $T_{\lambda}$.

\begin{Def}
\label{defc}
Let $\mathfrak{t}\subset\mathfrak{g}$ be a torus, centralized by a homomorphism $\phi:\Gamma\to G$. Among the roots of the adjoint action of $\mathfrak{t}$ on $\mathfrak{g}$, let $P$ be the subset of pure imaginary roots $\lambda$ such that $\rho_{\lambda}$ is a maximal representation with $T_{\lambda}>0$. Say $\mathfrak{t}$ is \emph{balanced} with respect to $\phi$ if 0 belongs to the interior of the sum of the convex hull of the imaginary parts of elements of $P$ and the linear span of the real and imaginary parts of roots not in $\pm P$.
\end{Def}

Here is a necessary and sufficient condition for flexibility, for surface groups of sufficiently large genus.

\begin{sthm}
\label{flexssimple}
Let $G$ be a semisimple real algebraic group. Let $\Gamma$ be the fundamental group of a closed surface of genus $\geq 2\mathrm{dim}(G)^2$. Let $\phi:\Gamma\to G$ be a homomorphism with reductive Zariski closure. Then $\phi$ is flexible if and only if $\mathfrak{c}$, the center of the centralizer of $\phi(\Gamma)$, is balanced with respect to $\phi$.
\end{sthm}

Using a result by M. Burger, A. Iozzi and A. Wienhard leads to the following consequence.

\begin{scor}
\label{rigidtt}
Let $G$ be a reductive real algebraic group. Let $\Gamma$ be the fundamental group of a closed surface of genus $\geq 2\mathrm{dim}(G)^2$. Assume there exists a non flexible homomorphism $\phi:\Gamma\to G$. Then the Zariski closure of $\phi(\Gamma)$ admits a transitive isometric action on a tube type Hermitian symmetric space, and the action of $\Gamma$ on this space is a maximal representation.
\end{scor}

It is not that easy to decide whether the balance condition is satisfied or not by a given homomorphism $\phi:\Gamma\to G$. For rank $1$ simple groups, it is possible, and a complete answer to Question \ref{question1} for large genus can be given\footnote{An other class, classical simple groups
will be treated in the companion paper \cite{classical}}. Up to the genus restriction, this is a converse to Toledo's rigidity theorem \ref{Toledo}.

\begin{scor}
\label{rankone}
Let $G$ be a rank one almost simple Lie group. Let $\Gamma$ be the fundamental group of a closed surface of genus $\geq 2\mathrm{dim}(G)^{2}$. Let $\phi:\Gamma\to G$ be a homomorphism. Then $\phi$ is flexible, unless $G$ is $SU(m,1)$ and $\phi(\Gamma)$ is discrete, cocompact in a conjugate of $S(U(1,1)\times U(m-1))\subset SU(m,1)$.
\end{scor}

\begin{rem}
\label{remgenus}
It is likely that the restriction on genus is irrelevant. Nevertheless, our arguments dwelve heavily on it.
\end{rem}

\begin{rem}
\label{remhiggs}
Higgs bundles.
\end{rem}
Character varieties of surface groups can be described in terms of Higgs bundles, leading to a purely algebro-geometric approach to Question \ref{question1}. It gives detailed informations on the topology of character varieties (see \cite{Hit1}, \cite{Hit2},\cite{BGG}) in many cases. But, as far as we can see, the published results do not seem to contain a complete anwer to Question \ref{question1}.

\subsection{Tools}

Our main tools were already known and used by W. Goldman in his 1985 paper, \cite{G1}. He was able to reduce the local study of the representation variety $Hom(\Gamma,G)$ to cohomology calculations : near the conjugacy class of a reductive homomorphism $\phi$, $Hom(\Gamma,G)$ identifies with its quadratic approximation, the set of cocycles in $Z^1 (\Gamma,\mathfrak{g}_{Ad\circ\phi})$ whose cup-square vanishes. Then, in an example, he was able to compute the cup-square map, thanks to an avatar of the index theorem due to W. Meyer, \cite{Meyer}.

\medskip

The first step in our investigation of flexibility deals with centers of centralizers. They split the Lie algebra $\mathfrak{g}$ into symplectic root spaces $\mathfrak{g}_{\lambda,\br}$, the cup-square map splits accordingly. On each $H^1 (\Gamma,\mathfrak{g}_{\lambda,\br})$, the cup-square map is {\em scalar}, i.e. of the form $Q_{\lambda}t_{\lambda}$ where $t_{\lambda}\in H^2 (\Gamma,\mathfrak{g}_{\lambda,\br})$ is a vector, and $Q_{\lambda}$ a scalar valued quadratic form whose signature can be interpreted as a Toledo invariant. The quadratic approximation has a dense set of smooth points (and flexibility holds) unless sufficiently many of the $Q_{\lambda}$'s are definite, i.e. the representation of $\Gamma$ on $\mathfrak{g}_{\lambda,\br}$ is maximal. Tube-type Hermitian domains then arise from Burger-Iozzi-Wienhard's results from \cite{BIW}. It turns out that balancedness is necessary and sufficient simultaneously for existence and density of smooth points in the quadratic approximation, this is the origin of Theorem \ref{thm1}.

The second step is a direct construction of cohomology classes which
shows that representations with semisimple centralizers are
flexible. This is based on some linear algebra of symplectic vector
spaces, and involves a characteristic class computation. More generally, one shows that the presence of a Levi factor in the centralizer does not affect the issue of flexibility.

Further steps allow reduction from general representations to reductive ones, and from reductive ambient groups to semisimple ones.

\subsection{Organization of the paper}

Section \ref{fuchsian} explains the dimension count that reduces density of Zariski density to smoothness of the character variety. The short section \ref{GMS} recalls W. Goldman's approach to rigidity. In section \ref{toral}, the role of the center of the centralizer is unveiled. Levi factors of centralizers appear in section \ref{levi}. The flexibility of the trivial representation (Theorem \ref{thm0}) appears in subsection \ref{trivialh}. The flexibility statement in Theorem \ref{flexssimple} is proven in subsection \ref{proofflexsimpleflex}. Due to complications caused by finite centers, the converse (rigidity statement) requires lifting homomorphisms in central extensions (section \ref{ext}). It also requires flexibility of representations with abelian, or, more generally, amenable image (section \ref{amen}), so the full proof of Theorem \ref{flexssimple} is postponed until subsection \ref{proofflexsimplerig}. The last reduction, to cover non reductive homomorphisms, comes in section \ref{nonreductive}, where the proofs of Theorem \ref{thm1} and Corollary \ref{rigidtt} are completed. Section \ref{one} deals with the rank one case and Corollary \ref{rankone}. Two appendices collect results not easily accessible in the litterature, W. Meyer's signature formula and the structure of centralizers of reductive subgroups.

The authors thank Yves Benoist, Jean-Louis Clerc, Patrick Eberlein, Philippe Gilles, for their help in dealing with algebraic groups.

\section{Dimension considerations}
\label{fuchsian}

The point of this section is to prove that Zariski dense homomorphisms $\Gamma\to G$, $\Gamma$ the fundamental group of a closed surface of genus $>1$, $G$ the group of real points of a reductive real algebraic group, form a dense subset of $\mathcal{S}$ when
the genus is high enough.

\subsection{The virtual dimension \texorpdfstring{$\mathrm{vdim}(Hom(\Gamma,G))$}{}}
\label{thevirtual}

Let $\Gamma$ a finitely generated group. Fix once and for all a finite generating set $S\subset\Gamma$. Let $\mathbb{F}$ denote the free group generated by $S$, and $N\subset\mathbb{F}$ the set of relations, i.e. the kernel of the tautological homomorphism $\mathbb{F}\to \Gamma$.
Let $G$ be a Lie group. Let $F:G^{S}\to G^{N}$ map $\phi\in G^{S}$, viewed as a homomorphism $\mathbb{F}\to G$, to the collection of all $\phi(n)$, $n\in N$. Then $Hom(\Gamma,G)=F^{-1}(\mathbf{e})$, where $\mathbf{e}$ is the constant collection $(e,e,\ldots)$. This provides $Hom(\Gamma,G)$ with the structure of a real analytic set.

At a point $\phi$, the Zariski tangent space of $Hom(\Gamma,G)$ identifies with the space $Z^{1}(\Gamma,\mathfrak{g}_{ad\circ\phi})$ of $1$-cocycles with values in the Lie algebra $\mathfrak{g}$. The dimension of the Zariski tangent space is
\begin{eqnarray*}
\mathrm{dim}(Z^{1}(\Gamma,\mathfrak{g}_{ad\circ\phi}))
&=&\mathrm{dim}(H^{1}(\Gamma,\mathfrak{g}_{ad\circ\phi}))+\mathrm{dim}(B^{1}(\Gamma,\mathfrak{g}_{ad\circ\phi}))\\
&=&\mathrm{dim}(H^{1}(\Gamma,\mathfrak{g}_{ad\circ\phi}))+\mathrm{dim}(C^{0}(\Gamma,\mathfrak{g}_{ad\circ\phi}))\\
&&-\mathrm{dim}(H^{0}(\Gamma,\mathfrak{g}_{ad\circ\phi})).
\end{eqnarray*}
For $\Gamma$ a closed surface group, and $\rho$ a linear representation of $\Gamma$ on a vector space $V_\rho$, elementary homological algebra gives
\begin{eqnarray*}
\mathrm{dim}(H^0 (\Gamma,V_\rho))-\mathrm{dim}(H^1 (\Gamma,V_\rho))+\mathrm{dim}(H^2 (\Gamma,V_\rho))
=\chi(\Gamma)\mathrm{dim}(V_\rho).
\end{eqnarray*}
For $\rho=ad\circ\phi$, this yields
\begin{eqnarray*}
\mathrm{dim}(Z^{1}(\Gamma,\mathfrak{g}_{ad\circ\phi}))&=&\mathrm{dim}(H^{2}(\Gamma,\mathfrak{g}_{ad\circ\phi}))+(1-\chi(\Gamma))\mathrm{dim}(G).
\end{eqnarray*}

If $G$ is connected reductive with radical $R$ (which coincides with its center), pick a semisimple Levi factor $S$. Then $\mathfrak{g}=\mathfrak{s}\oplus\mathfrak{r}$, and
\begin{eqnarray*}
H^2 (\Gamma,\mathfrak{g}_{ad\circ\phi})=H^2 (\Gamma,\mathfrak{s}_{ad\circ\phi})\oplus H^2 (\Gamma,\br)\otimes\mathfrak{r}.
\end{eqnarray*}
By semi-continuity, if $H^2 (\Gamma,\mathfrak{s}_{ad\circ\phi})=0$, this still holds for neighboring homomorphisms. Therefore the rank of the differential of $F$ is locally constant, and $Hom(\Gamma,G)$ is a smooth analytic manifold in a neighborhood of $\phi$, of dimension $\mathrm{dim}(\mathfrak{r})+(1-\chi(\Gamma))\mathrm{dim}(G)$. This suggests the following definition.

\begin{Def}
\label{defvirtual}
Let $G$ be a reductive real algebraic group with radical $R$. Let $\Gamma$ be a closed surface group. Define the {\em virtual dimension} of $Hom(\Gamma,G)$ as
\begin{eqnarray*}
\mathrm{vdim}(Hom(\Gamma,G))=(1-\chi(\Gamma))\mathrm{dim}(G)+\mathrm{dim}(R).
\end{eqnarray*}
Say a homomorphism $\phi$ is a {\em smooth} point of $Hom(\Gamma,G)$ if
$$\mathrm{dim}(Z^1 (\Gamma,\mathfrak{g}_{ad\circ\phi}))=\mathrm{vdim}(Hom(\Gamma,G)).$$
Let $\mathcal{S}$ denote the set of smooth points of $Hom(\Gamma,G)$.
\end{Def}
Since smoothness is defined in terms of the Zariski tangent space only, it is invariant under analytic isomorphisms.

\subsection{Zariski dense homomorphisms belong to \texorpdfstring{$\mathcal{S}$}{}}
\label{virtual}

Let $\Gamma$ be a surface group and $\rho$ a linear representation of $\Gamma$. Poincar\'e duality gives an isomorphism of $H^2 (\Gamma,V_\rho)$ to $H^0 (\Gamma,V_\rho^*)$, where $\rho^*$ denotes the contragredient representation. Since $H^0 (\Gamma,V_\rho)$ counts $\Gamma$-invariant vectors in $V$,  $H^2 (\Gamma,V_\rho)$ counts $\Gamma$-invariant vectors in $V^*$.

If $\phi:\Gamma\to G$ is a homomorphism to a semisimple group $G$, invariant vectors in the adjoint representation $\rho=Ad\circ\phi$ form the Lie algebra of the centralizer $Z_G (\phi(\Gamma))$ of $\phi(\Gamma)$ in $G$. Since the adjoint representation preserves a nondegenerate quadratic form, the Killing form, it is isomorphic to its contragredient. Therefore $H^2
(\Gamma,\mathfrak{g}_{Ad\circ\phi})$ is isomorphic to $H^0
(\Gamma,\mathfrak{g}_{Ad\circ\phi})$.

\begin{lemma}
\label{nocentralizer}
Let $G$ be semisimple, let $\Gamma$ be a surface group. Let $\phi:\Gamma\ra G$ be a homomorphism such that $\phi(\Gamma)$ has a discrete centralizer. Then $\phi$ belongs to the set $\mathcal{S}$ of smooth homomorphisms (see Definitions \ref{defrigid} and \ref{defvirtual}). In particular, Zariski dense homomorphisms belong to $\mathcal{S}$.
\end{lemma}

\begin{pf}
If $\phi(\Gamma)$ has a discrete centralizer, then $H^2
(\Gamma,\mathfrak{g}_{Ad\circ\phi})=H^0
(\Gamma,\mathfrak{g}_{Ad\circ\phi})$ $=0$, which implies that $\phi$ is a smooth point of $Hom(\Gamma,G)$.
\end{pf}

When $G$ is merely reductive, there is an analogue of Lemma \ref{nocentralizer}, Proposition \ref{centralizercenter}. It relies on properties of central extensions which are postponed until section \ref{ext}.

\subsection{A crude upper bound on \texorpdfstring{$\mathrm{dim}(Hom(\Gamma,G))$}{}}

Let $\Gamma$ be a surface group and $\rho$ a linear representation of $\Gamma$. Then $\mathrm{dim}(H^0 (\Gamma,V_\rho))\leq\mathrm{dim}(V_{\rho})$, $\mathrm{dim}(H^0 (\Gamma,V_\rho))\leq\mathrm{dim}(V_{\rho}^{*})$, thus
\begin{eqnarray*}
\mathrm{dim}(H^1 (\Gamma,V_\rho)) \leq (2-\chi(\Gamma))\mathrm{dim}(V_{\rho}).
\end{eqnarray*}
and
\begin{eqnarray}
\label{inecrude}
\mathrm{dim}(Z^1 (\Gamma,V_\rho)) \leq (3-\chi(\Gamma))\mathrm{dim}(V_{\rho}).
\end{eqnarray}

\begin{lemma}
\label{crude}
Let $G$ be a Lie group. Let $\Gamma$ be a finitely generated group. Then, in a neighborhood of the conjugacy class of a homomorphism $\phi:\Gamma\to G$,
\begin{eqnarray*}
\mathrm{dim}(Hom(\Gamma,G))\leq(3-\chi(\Gamma))\mathrm{dim}(G).
\end{eqnarray*}
\end{lemma}

\begin{pf}
Since one can locally embed $Hom(\Gamma,G)$ in its Zariski tangent space, its dimension is at most $\mathrm{dim}(Z^{1}(\Gamma,\mathfrak{g}_{ad\circ\phi}))\leq (3-\chi(\Gamma))\mathrm{dim}(G)$, by inequality (\ref{inecrude}).
\end{pf}

\subsection{Zariski dense homomorphisms are dense in \texorpdfstring{$\mathcal{S}$}{}}

\begin{co}
\label{dense}
Let $G$ be a connected reductive real algebraic group of dimension $>2$. Let $\Gamma$ be the fundamental group of a closed surface of genus $g$. Assume that $g\geq \mathrm{dim}(G)^2$. Then the set of non Zariski dense homomorphisms $\Gamma\ra G$ has dimension $<-\chi(\Gamma)\mathrm{dim}(G)$. In particular, Zariski dense representations form a dense subset of $\mathcal{S}$.
\end{co}

\begin{pf}
Connected proper subgroups of $G$ are determined by their Lie algebras, which are linear subspaces in $\mathfrak{g}$. Algebraic subgroups of $G$ have finitely many connected components. Therefore proper algebraic subgroups of $G$ come in countably many families, each of which has dimension less than $\mathrm{dim}(G)^2$. The set of pairs $(H,\phi)$ where $H\subset G$ is a proper algebraic subgroup and $\phi:\Gamma\to H$ is a homomorphism has dimension at most $\mathrm{dim}(G)^2+(3-\chi(\Gamma))(\mathrm{dim}(G)-1)$. If $g\geq \mathrm{dim}(G)^2$,
\begin{eqnarray*}
&&-\chi(\Gamma)\mathrm{dim}(G)-(\mathrm{dim}(G)^2+(3-\chi(\Gamma))(\mathrm{dim}(G)-1))\\
&=&2g-2-\mathrm{dim}(G)^2-3\mathrm{dim}(G)+3\\
&\geq&\mathrm{dim}(G)^2-3\mathrm{dim}(G)+1>0.
\end{eqnarray*}
Since $\mathcal{S}$ is a smooth $\mathrm{vdim}(Hom(\Gamma,G))$-dimensional manifold, and
$$\mathrm{vdim}(Hom(\Gamma,G))\geq-\chi(\Gamma)\mathrm{dim}(G),$$
a subset of smaller dimension has empty interior in it.
\end{pf}

\begin{rem}
\label{remconnected}
The connectedness assumption in Corollary \ref{dense} is necessary.
\end{rem}
Indeed, if $G'\subset G$ is a proper open subgroup, and genus is large enough, homomorphisms $\Gamma\to G'$ form a proper open and closed subset of $Hom(\Gamma,G)$, there are smooth ones but none of them is Zariski dense.

This motivates the following generalization of Corollary \ref{dense}.

\begin{Prop}
\label{dense'}
Let $G$ be a reductive real algebraic group of dimension $>2$. Let $\Gamma$ be the fundamental group of a closed surface of genus $g$. Assume that $g\geq \mathrm{dim}(G)^2$. Let $\phi:\Gamma\to G$ be a smooth homomorphism. Let $G'$ denote the smallest open subgroup of $G$ containing $\phi(\Gamma)$. Then homomorphisms whose Zariski closure is $G'$ are dense in a neighborhood of $\phi$ in $Hom(\Gamma,G)$.
\end{Prop}

\begin{pf}
$G'$ is well defined since $G$ has only finitely many connected components. As $Hom(\Gamma,G')$ is open in $Hom(\Gamma,G)$, $\phi$ is smooth as a homomorphism $\Gamma\to G'$, so we can assume that $G=G'$. There are finitely many proper open subgroups $G''$ in $G$. Each $Hom(\Gamma,G'')$ is closed in $Hom(\Gamma,G)$, so is their union. Thus a neighborhood $\mathcal{V}$ of $\phi$ in $Hom(\Gamma,G)$ consists of homomorphisms which are not contained in any proper open subgroup of $G$. As observed in the proof of Corollary \ref{dense}, in $\mathcal{V}$, homomorphisms which are contained in proper closed subgroups of empty interior form a set of lower dimension. Thus Zariski dense homomorphisms are dense in $\mathcal{V}$.
\end{pf}

\begin{ex}
\label{exirreducible}
Let $\Gamma$ be the fundamental group of a closed surface of genus $g$. Let $\phi:\Gamma\to Sl(n,\br)$ be an irreducible representation. Then $\phi$ is smooth and thus flexible. If furthermore $g\geq n^4$, $\phi$ can be deformed to become Zariski dense.
\end{ex}

\begin{pf}
An element of the centralizer of $\phi(\Gamma)$ is an invertible linear map which intertwines the representation $\phi$ with itself. According to Schur's lemma, it has to be proportional to identity. In other words, the centralizer of $\phi(\Gamma)$ in $Sl(n,\br)$ is finite. Thus $\phi$ belongs to the open set $\mathcal{S}$. In particular, $\phi$ is flexible. If $g\geq n^4$, Corollary \ref{dense} applies, and Zariski dense homomorphisms are dense in a neighborhood of $\phi$.
\end{pf}

\section{Second order calculations}
\label{GMS}

The second order integrability condition for infinitesimal
deformations at $\phi$ of representations of a group $\Gamma$ in a
Lie group $G$ can be expressed in terms of the {\em cup-product}, a
symmetric bilinear map
\begin{eqnarray*}
[\cdot \smile\cdot]:H^1(\Gamma, \mathfrak{g}_{Ad\circ\phi})\to H^2(\Gamma,
\mathfrak{g}_{Ad\circ\phi}).
\end{eqnarray*}
For $u\in Z^1(\Gamma, \mathfrak{g}_{Ad\circ\phi})$,
$$
[u\smile u](\alpha,\beta)=[u(\alpha),Ad\circ\phi(\alpha)u(\beta)].
$$
It is well-known, \cite{NR}, that for a representation $\phi$ from
$\Gamma$ to a reductive group $G$, if there exists a smooth path
$\phi_t$ in $Hom(\Gamma, G)$ which is tangent to $u\in Z^1(\Gamma,
\mathfrak{g}_{Ad\circ\phi})$, then $[u\smile u]=0$. According to W. Goldman (Theorem 3 in \cite{G1}, \cite{GM}), for surface groups, this necessary condition is also
sufficient. Here is a more general statement, due to W. Goldman, J. Millson and C. Simpson, following an idea of P. Deligne.

\begin{Def}
Let $\Gamma$ be a finitely generated group, let $\phi:\Gamma\to G$ be a homomorphism to a reductive Lie group $G$. The \emph{quadratic model} for the representation space near $\phi$ is
\begin{eqnarray*}
\mathcal{Q}_{\phi}=\{u\in Z^{1}(\Gamma, \mathfrak{g}_{Ad\circ\phi})\,|\,[u\smile u]=0~\mathrm{in}~H^{2}(\Gamma, \mathfrak{g}_{Ad\circ\phi})\}.
\end{eqnarray*}
\end{Def}

\begin{thm}
\label{dgms}
{\em (W. Goldman and J. Millson, \cite{GM3}, Theorem 1, for bounded homomorphisms, C. Simpson, \cite{S}, corollary 2.4, for general reductive homomorphisms).} Let $G$ be a semi-simple Lie group. Let $\Gamma$ be the fundamental group of a compact K\"ahler manifold. If the Zariski closure of $\phi(\Gamma)$ is reductive, the analytic germ of the algebraic variety $Hom(\Gamma,G)$ at $\phi$ is equivalent to the analytic germ of $\mathcal{Q}_{\phi}$ at $0$.
\end{thm}

Under the analytic isomorphism of Theorem \ref{dgms}, smooth points of $Hom(\Gamma,G)$ are mapped to smooth points of the quadratic model, i.e. points $v\in \mathcal{Q}_{\phi}$ where the map $u\mapsto [u\smile u]$, $Z^{1}(\Gamma, \mathfrak{g}_{Ad\circ\phi})\ra H^{2}(\Gamma, \mathfrak{g}_{Ad\circ\phi})$, is a submersion.

\begin{rem}
Theorem \ref{dgms} implies that locally, $Hom(\Gamma,G)$ can be defined by homogeneous quadratic equations. This is not sufficient to imply Theorem \ref{thm1}.
\end{rem}
Indeed, consider the following system of $2$ homogeneous quadratic equations on $\br^4$ :
\begin{eqnarray*}
\left\{\begin{array}{ccc}
x^2 +y^2 -tz-t^2 &=&0, \\
x^2 +y^2+z^2 -t^2&=&0.
\end{array}\right.
\end{eqnarray*}
The solution set $\mathcal{Q}$ is the union of the line $L=\{x=0,\,y=0,\,z=-t\}$ (which is $1$-dimensional) and of the quadratic cone $C=\{x^2 +y^2 -t^2 =0,\,z=0\}$ (which is $2$-dimensional). Using the terminology of Definitions \ref{defrigid}, \ref{defvirtual}, the smooth points of $\mathcal{Q}$, as well as the flexible ones, are the points of $C\setminus\{0\}$, the rigid points of $\mathcal{Q}$ are the points of $L\setminus\{0\}$, and $0$ is neither flexible nor rigid.

\section{Toral centralizers}
\label{toral}

\subsection{Toledo invariants attached to a torus}

\begin{Prop}
\label{torus}
Let $G$ be a semisimple real algebraic group. Let $\mathfrak{t}$ be a (non necessarily maximal) torus in $\mathfrak{g}$. After complexification, $\mathfrak{g}$ splits as
\begin{eqnarray*}
\mathfrak{g}\otimes\bc=\mathfrak{g}_{0}\oplus\bigoplus_{\lambda\in\Lambda}\mathfrak{g}_{\lambda}.
\end{eqnarray*}
For a nonzero root $\lambda\in\Lambda$, the real vector space
\begin{eqnarray*}
\mathfrak{g}_{\lambda,\br}=\mathfrak{g}\cap((\mathfrak{g}_{\lambda}\oplus\mathfrak{g}_{-\lambda})+(\mathfrak{g}_{\bar{\lambda}}\oplus\mathfrak{g}_{-\bar{\lambda}})).
\end{eqnarray*}
carries a natural complex valued alternating 2-form $\Omega_{\lambda}$. If $t_{\lambda}\in\mathfrak{t}\otimes\bc$ denotes the complex root vector defined by $\lambda(t)=t_{\lambda}\cdot t$ (complexified Killing inner product), then, for $X$, $Y\in \mathfrak{g}_{\lambda,\br}$,
\begin{eqnarray*}
[X,Y]^{\mathfrak{t}}=\Re e(\Omega_{\lambda}(X,Y)t_{\lambda}),
\end{eqnarray*}
(Killing orthogonal projection on $\mathfrak{t}$).
Also, $\Omega_{-\lambda}=-\Omega_{\lambda}$, $\Omega_{\bar{\lambda}}=\overline{\Omega_{\lambda}}$.

If $\lambda\not=0$ is real, so are $t_{\lambda}$ and $\Omega_{\lambda}$. $\Omega_{\lambda}$ is a symplectic structure on $\mathfrak{g}_{\lambda,\br}$.

If $\lambda\not=0$ is pure imaginary, so are $t_{\lambda}$ and $\Omega_{\lambda}$. The imaginary part $\Im m(\Omega_{\lambda})$ is a symplectic structure on $\mathfrak{g}_{\lambda,\br}$.

If $\lambda$ is neither real nor pure imaginary, $\mathfrak{g}_{\lambda,\br}$ carries a natural complex structure $J_{\lambda}$ and $\Omega_{\lambda}$ is $\bc$-bilinear with respect to $J_{\lambda}$. It defines a complex symplectic structure on $(\mathfrak{g}_{\lambda,\br},J_{\lambda})$.
\end{Prop}

\begin{pf}
By assumption, after complexification, $\mathfrak{g}$ splits into root spaces
\begin{eqnarray*}
\mathfrak{g}\otimes\bc=\mathfrak{g}_{0}\oplus\bigoplus_{\lambda\in\Lambda}\mathfrak{g}_{\lambda},
\end{eqnarray*}
where $\Lambda\subset\mathfrak{t}^* \otimes\bc$ is a finite set of roots, and $\mathfrak{g}_{\lambda}=\{X\in\mathfrak{g}\otimes\bc\,|\,\forall t\in \mathfrak{t},\,ad_t (X)=\lambda(t)X\}$. For $t\in\mathfrak{t}$, consider the alternating form on $\mathfrak{g}$ defined by
\begin{eqnarray*}
\Omega_t (X,Y)=t\cdot[X,Y].
\end{eqnarray*}
Since the Killing form is bi-invariant, for any $X$, $Y\in\mathfrak{g}$, $ad_X (Y)\cdot t +Y\cdot ad_X (t)=0$, thus
\begin{eqnarray*}
\Omega_t (X,Y)=ad_t (X)\cdot Y
\end{eqnarray*}
for all $t\in\mathfrak{t}$.  Since
$[\mathfrak{g}_{\alpha},\mathfrak{g}_{\beta}]\subset
\mathfrak{g}_{\alpha+\beta}$, the subspaces $\mathfrak{g}_{\alpha}$
and $\mathfrak{g}_{\beta}$ are orthogonal with respect to the
Killing form and to $\Omega_t$ unless $\alpha+\beta=0$. Since
$ad_t$'s are skewsymmetric, all nonzero roots occur in pairs
$(\lambda,-\lambda)$.  For nonzero $\lambda\in\Lambda$, set
$\mathfrak{g}_{\pm\lambda}=\mathfrak{g}_{\lambda}\oplus\mathfrak{g}_{-\lambda}$.
Distinct $\mathfrak{g}_{\pm\lambda}$'s are Killing and
$\Omega_t$-orthogonal. More precisely let $X\in
\mathfrak{g}_\alpha,Y\in \mathfrak{g}_\beta$ so that
$\alpha+\beta\neq 0$ and $Z\in \mathfrak{z}$ so that
$(\alpha+\beta)(Z)\neq 0$. From $[X,Y]\cdot Z=-\beta(Z)X\cdot
Y=\alpha(Z)X\cdot Y$ we get $(\alpha+\beta)(Z) X\cdot Y=0$ for any
$X,Y$. So $\mathfrak{g}_\alpha$ and $\mathfrak{g}_\beta$ are
orthogonal. Similar for $\Omega_Z$. It follows that the Killing form
restricted to each $\mathfrak{g}_{\pm\lambda}$ and to
$\mathfrak{g}_{0}$ is nondegenerate. Furthermore, $\Omega_t$
restricted to each $\mathfrak{g}_{\pm\lambda}$ such that
$\lambda(t)\not=0$ is nondegenerate, and the nullspace of $\Omega_t$
equals the sum of $\mathfrak{g}_{\lambda}$ such that $\lambda(t)=0$.

Let $\lambda$ be a nonzero root. Let
\begin{eqnarray*}
\mathfrak{g}_{\lambda,\br}:=\mathfrak{g}\cap(\mathfrak{g}_{\pm\lambda}+\mathfrak{g}_{\pm\bar{\lambda}})=\mathfrak{g}\cap(\mathfrak{g}_{\lambda}+\mathfrak{g}_{-\lambda}+\mathfrak{g}_{\bar{\lambda}}+\mathfrak{g}_{\bar{\lambda}}),
\end{eqnarray*}
(the sum is direct only if $\lambda$ is neither real nor pure imaginary). Let $X=X_{\lambda}+X_{-\lambda}$, $Y=Y_{\lambda}+Y_{-\lambda}\in\mathfrak{g}_{\pm\lambda}$. For $t\in\mathfrak{t}$,
\begin{eqnarray*}
ad_{t}(X)=\lambda(t)(X_{\lambda}-X_{-\lambda}),
\end{eqnarray*}
\begin{eqnarray*}
\Omega_t (X,Y)=ad_{t}(X)\cdot Y=\lambda(t)(X_{\lambda}-X_{-\lambda})\cdot(Y_{\lambda}+Y_{-\lambda})
\end{eqnarray*}
is a complex multiple of the complex symplectic structure
\begin{eqnarray*}
\omega_{\lambda} (X,Y)=(X_{\lambda}-X_{-\lambda})\cdot(Y_{\lambda}+Y_{-\lambda}).
\end{eqnarray*}

If $\lambda$ is neither real nor pure imaginary,  the splitting
$\mathfrak{g}_{\lambda,\br}=\mathfrak{g}_{\pm\lambda}\oplus\overline{\mathfrak{g}_{\pm\lambda}}$
defines a complex structure $J_{\lambda}$ on
$\mathfrak{g}_{\lambda,\br}$ (its complexification equals
multiplication by $ i$ on $\mathfrak{g}_{\pm\lambda}$ and by $-i$ on
$\mathfrak{g}_{\pm\bar\lambda}$). Then $\omega_{\lambda}$ defines a
$J_{\lambda}$-bilinear complex valued alternating 2-form
$\Omega_{\lambda}$ on $\mathfrak{g}_{\lambda,\br}$ as follows :
\begin{eqnarray*}
\Omega_{\lambda}(X+\bar{X},Y+\bar{Y})=2\omega_{\lambda} (X,Y).
\end{eqnarray*}
Let $t_{\lambda}\in\mathfrak{t}\otimes\bc$ be the vector defined by
\begin{eqnarray*}
t_{\lambda}\cdot t=\lambda(t)\quad\textrm{ for all }t\in\mathfrak{t}\otimes\bc.
\end{eqnarray*}
Then, for $X$, $Y\in\mathfrak{g}_{\pm\lambda}$,
\begin{eqnarray*}
t\cdot[X,Y]=\Omega_t (X,Y)=\lambda(t)\omega_{\lambda} (X,Y)=t\cdot(\omega_{\lambda} (X,Y)t_{\lambda}),
\end{eqnarray*}
thus
\begin{eqnarray*}
[X,Y]^{\mathfrak{t}\otimes\bc}=\omega_{\lambda} (X,Y)t_{\lambda},
\end{eqnarray*}
and, for real $X+\bar{X}$, $Y+\bar{Y}\in\mathfrak{g}_{\lambda,\br}$,
since $$ t\cdot [X+\bar{X},Y+\bar{Y}]=t\cdot[X,Y]+t\cdot[\bar X,\bar
Y]=\lambda(t)\omega_\lambda(X,Y)+\bar\lambda(t)\omega_{\bar\lambda}(\bar
X,\bar Y)$$
\begin{eqnarray*}
 [X+\bar{X},Y+\bar{Y}]^{\mathfrak{t}}
&=&([X,Y]+[\bar{X},\bar{Y}])^{\mathfrak{t}}\\
&=&2\Re e(\omega_{\lambda} (X,Y)t_{\lambda})\\
&=&\Re e(\Omega_{\lambda} (X+\bar{X},Y+\bar{Y})t_{\lambda}).
\end{eqnarray*}

If $\lambda\in\Lambda$ is real, $\Omega_{\lambda}$ is real and symplectic.

If $\lambda\in\Lambda$ is pure imaginary, since $t\cdot [X+\bar
X,Y+\bar Y]=\lambda(t)\Omega_\lambda(X+\bar X,Y+\bar Y)$ is real,
$\Omega_{\lambda}$ is pure imaginary and $\Im m(\Omega_{\lambda})$
is symplectic where
$$\Omega_\lambda(X+\bar X,Y+\bar Y)=(X-\bar X)\cdot(Y+\bar Y)$$ is
pure imaginary.
\end{pf}

\begin{Prop}
\label{abel} Let $G$ be a semisimple Lie group. Let $\Gamma$ be a
surface group.  Let $\phi:\Gamma\to G$ be a homomorphism whose
Zariski closure is reductive. Let $\mathfrak{c}$ denote the center
of its centralizer. Among the roots of the adjoint action of
$\mathfrak{c}$ on $\mathfrak{g}$, pick a subset
$\Lambda'\subset\Lambda$ containing exactly one representative of
each set $\{\lambda,-\lambda,\bar{\lambda},-\bar{\lambda}\}$. Pick a
$G^{\mathfrak{c}}$-invariant splitting
$\mathfrak{z}=\mathfrak{c}\oplus\mathfrak{s}$ with $\mathfrak{s}$ a
semisimple Lie subalgebra. Since
$H^{2}(\Gamma,\mathfrak{g})=\mathfrak{z}$, the quadratic map
$Q:H^{1}(\Gamma,\mathfrak{g})\ra H^{2}(\Gamma,\mathfrak{g})$ splits
as $Q=Q^{\mathfrak{c}}+Q^{\mathfrak{s}}$. Then, on
\begin{eqnarray*}
H^{1}(\Gamma,\mathfrak{g})=H^{1}(\Gamma,\mathfrak{g}_{0})
\oplus\bigoplus_{\lambda\in\Lambda'}H^{1}(\Gamma,\mathfrak{g}_{\lambda,\br}),
\end{eqnarray*}
the $\mathfrak{c}$-valued quadratic map $Q^{\mathfrak{c}}$ vanishes
on $H^{1}(\Gamma,\mathfrak{g}_{0})$  and splits as a direct sum of
nondegerate scalar quadratic forms
\begin{eqnarray*}
Q^{\mathfrak{c}}=\bigoplus_{\lambda\in\Lambda'} \Re e(Q_{\lambda}t_{\lambda}).
\end{eqnarray*}

For real or pure imaginary $\lambda$, $\mathfrak{g}_{\lambda,\br}$
is a symplectic vector space, whence a symplectic representation
$\rho_{\lambda}$ of $\Gamma$, with Toledo invariant $T_{\lambda}$.
Then the real valued quadratic form $Q_{\lambda}$ (resp. $\Im
m(Q_{\lambda})$) lives on a vector space $H^1
(\Gamma,\mathfrak{g}_{\lambda,\br})$ whose dimension is equal to
$-\chi(\Gamma)\mathrm{dim}(\mathfrak{g}_{\lambda,\br})$. The
signature of $Q_{\lambda}$ is equal to $4T_{\lambda}$.

If $\lambda$ is real, $T_{\lambda}=0$ automatically.

If $\lambda$ is neither real nor pure imaginary, $\mathfrak{g}_{\lambda,\br}$ has a complex structure, $Q_{\lambda}$ is a $\bc$-bi-linear form.

\end{Prop}

\begin{pf}
Since $\phi$ is reductive, its centralizer is reductive (Corollary \ref{centralizerisreductive}), and so is its center $\mathfrak{c}$. Thus $\mathfrak{c}$ is a torus, and we can apply Proposition \ref{torus}.

The splitting of cohomology follows from the corresponding splitting
of $\mathfrak{g}$  and of the $\mathfrak{c}$-component of the Lie
bracket. For nonzero $\lambda$,
$H^{0}(\Gamma,\mathfrak{g}_{\lambda})=\{a\in
\mathfrak{g}_{\lambda}|Ad\circ\phi(\gamma)a=a\ \text{ for}\ \text {all}\
\gamma\in \Gamma\}$ and such $a$ must belong to $\mathfrak{g}_{0}$,
and since $\mathfrak{g}_{\lambda}\cap \mathfrak{g}_0 =0$,
$H^{0}(\Gamma,\mathfrak{g}_{\lambda})=0$. Hence
$\mathrm{dim}_{\bc}H^{1}(\Gamma,\mathfrak{g}_{\lambda})=-\chi(\Gamma)\mathrm{dim}_{\bc}(\mathfrak{g}_{\lambda}))>0$,
thus no term vanishes. On
$H^{1}(\Gamma,\mathfrak{g}_{\pm\lambda}+\mathfrak{g}_{\pm\bar\lambda})$,
set a complex valued
\begin{eqnarray*}
Q_{\lambda}=\Omega_{\lambda}(u\smile u).
\end{eqnarray*}
Then, for $u=u_0 +\sum_{\lambda\in\Lambda'}u_{\lambda}$, a
$\mathfrak{c}$-valued
\begin{eqnarray*}
Q^{\mathfrak{c}}(u)&=&[u\smile u]^{\mathfrak{c}}\\
&=&\sum_{\lambda\in\Lambda'} [u_{\lambda}\smile u_{\lambda}]^{\mathfrak{c}}\\
&=&\sum_{\lambda\in\Lambda'} \Re e(Q_{\lambda}t_{\lambda}).
\end{eqnarray*}
Note that $\mathfrak{g}_{\lambda}$ is $Ad\circ\phi$-invariant. Since for
$t\in \mathfrak{c},X\in \mathfrak{g}_{\lambda}$,
$$[t,Ad\circ\phi(\gamma)X]=[Ad\circ\phi(\gamma)t,Ad\circ\phi(\gamma)X]=Ad\circ\phi(\gamma)[t,X]=\lambda(t)(Ad\circ\phi(\gamma)X).$$

For real or pure imaginary $\lambda$, let $\rho_{\lambda}:\Gamma\ra Sp(\mathfrak{g}_{\lambda,\br},\Omega_{\lambda})$ denote the symplectic linear representation defined by $\phi$ composed with the adjoint representation restricted to $\mathfrak{g}_{\lambda,\br}$. Let $E_{\rho_{\lambda}}$ be the corresponding symplectic vectorbundle on $\Sigma$. From Meyer's signature formula (see Lemma \ref{index}), it follows that the signature of the quadratic form $Q_{\lambda}(u)=\int_{\Sigma}\Omega_{\lambda}(u\smile u)$ on $H^{1}(\Gamma,\mathfrak{g}_{\lambda,\br})$ is equal to $4c_1 (E_{\rho_{\lambda}})$. Finally, Lemma \ref{toledochern} equates the first Chern class $c_1 (E_{\lambda})$ to the Toledo invariant $T_{\lambda}$.

If $\lambda$ is real, $\rho_{\lambda}$ fixes the splitting $\mathfrak{g}_{\lambda,\br}=(\mathfrak{g}\cap\mathfrak{g}_{\lambda})\oplus(\mathfrak{g}\cap\mathfrak{g}_{-\lambda})$ as a sum of real Lagrangian subspaces. It follows from Lemma \ref{toledochern} that $T_{\lambda}=0$.

If $\lambda$ is neither real nor pure imaginary, the complex structure $J_{\lambda}$ on $\mathfrak{g}_{\lambda,\br}$ defines a complex structure on $H^{1}(\Gamma,\mathfrak{g}_{\lambda,\br})$ : $u\mapsto J_{\lambda}\circ u$, with respect to which the complex valued quadratic form $Q_{\lambda}$ is $\bc$-bi-linear.
\end{pf}

\begin{ex}
\label{excomplex}
Let $G$ be a complex Lie group. Then no root $\lambda$ is pure imaginary.
\end{ex}
Indeed, the centralizer of an arbitrary subset of $G$ is a complex Lie subgroup, its center $\mathfrak{c}$ is a complex vectorspace and roots are $\bc$-linear, they cannot be pure imaginary.

\subsection{The central quadratic model}

\begin{Def}
\label{defcentrqm}
Let $G$ be a semisimple Lie group. Let $\Gamma$ be a surface group. Let $\phi:\Gamma\to G$ be a homomorphism whose Zariski closure is reductive. Let $\mathfrak{c}$ denote the center of the centralizer of $\phi(\Gamma)$ in $\mathfrak{g}$. We call the space
\begin{eqnarray*}
\mathcal{Q}_{\phi}^{\mathfrak{c}}=\{u\in Z^{1}(\Gamma,\mathfrak{g})\,|\,[u\smile u]^{\mathfrak{c}}=0~in~H^{2}(\Gamma,\mathfrak{c})\}
\end{eqnarray*}
the \emph{central quadratic model} for $Hom(\Gamma,G)$ near $\phi$. Its {\em smooth points} are the points where $u\mapsto [u\smile u]^{\mathfrak{c}}$ is a submersion. The projection of $\mathcal{Q}_{\phi}^{\mathfrak{c}}$ in $H^{1}(\Gamma,\mathfrak{g})$ is {\em denoted by} $\bar{\mathcal{Q}}_{\phi}^{\mathfrak{c}}$.
\end{Def}

By definition, the central quadratic model contains the Goldman-Millson-Simpson quadratic model, and coincides with it if the centralizer of $\phi(\Gamma)$ is abelian. A smooth point of the Goldman-Millson-Simpson quadratic model is smooth as a point of the central quadratic model as well.

\begin{lemma}
\label{smoothcenter}
Let $G$ be a semisimple Lie group. Let $\Gamma$ be a surface group. Let $\phi:\Gamma\to G$ be a reductive homomorphism. Let $C$ denote the center of the centralizer of $\phi(\Gamma)$ in $G$. A real cohomology class $u=u_0 +\sum_{\lambda\in\Lambda'} u_{\lambda}\in H^{1}(\Gamma,\mathfrak{g})$ defines a smooth point of the central quadratic model if and only if the set of roots $\lambda$ such that $u_{\lambda}\not=0$ spans $\mathfrak{c}^{*}\otimes\bc$.
\end{lemma}

\begin{pf}
The differential of $Q^{\mathfrak{c}}$ at $u$ is
\begin{eqnarray*}
u'\mapsto \sum_{\lambda\in\Lambda'}\Re
e(d_{u_{\lambda}}Q_{\lambda}(u')t_{\lambda}).
\end{eqnarray*}
The image of this linear map is clearly contained in the linear span of the set of all real and imaginary parts of root vectors $t_{\lambda}$ such that $u_{\lambda}\not=0$. In fact, it is equal to it. Indeed, when $\lambda\not=0$ and $u_{\lambda}$ is non zero, the map
\begin{eqnarray*}
H^{1}(\Gamma,\mathfrak{g}_{\lambda,\br})\to H^{2}(\Gamma,\mathfrak{c}),\quad u'\mapsto\Re e(d_{u_{\lambda}}Q_{\lambda}(u')t_{\lambda})
\end{eqnarray*}
is surjective onto the linear span of $\Re e(t_{\lambda})$ and $\Im m(t_{\lambda})$. If $\lambda$ is real or pure imaginary, this is because $Q_{\lambda}$ (resp. $\Im m(Q_{\lambda})$) is nondegenerate as a real bilinear form. If $\lambda$ is neither real nor pure imaginary, $H^{1}(\Gamma,\mathfrak{g}_{\lambda,\br})$ has a complex structure and $Q_{\lambda}$ is nondegenerate as a $\bc$-bilinear form. We conclude that $d_{u}Q^{\mathfrak{c}}$ is surjective if and only if
\begin{eqnarray*}
\mathfrak{c}=\mathrm{span}_{\br}(\bigcup_{\lambda\not=0,\,u_{\lambda}\not=0}\{\Re e(t_{\lambda}),\Im m(t_{\lambda})\}).
\end{eqnarray*}
An equivalent formulation is
\begin{eqnarray*}
\mathfrak{c}^{*}\otimes\bc=\mathrm{span}_{\bc}(\{\lambda\,|\,u_{\lambda}\not=0\}).
\end{eqnarray*}
\end{pf}

\begin{rem}
Note that under this condition, not only is $\bar{\mathcal{Q}}_{\phi}^{\mathfrak{c}}$ smooth near $u$, but its quotient under the stabilizer of $u$ in $C$, which is a local model for $\chigg=Hom(\Gamma/G)/G$, is a smooth orbifold of the same dimension.
\end{rem}
Indeed, the action of $C$ is locally free near $u$. For each root $\lambda$ such that $u_{\lambda}\not=0$, at
least one of the (2 or 4) components (in $\mathfrak{g}_{\pm
\lambda}$ and $\mathfrak{g}_{\pm \bar\lambda}$) of $u_{\lambda}$ on
root spaces does not vanish. Let us denote it by $x_{\lambda}$. If
$t\in \mathfrak{c}$ satisfies $ad_t (u)=0$, then $ad_t
(x_{\lambda})=0$. This implies that one of $\lambda(t)$,
$-\lambda(t)$, $\bar{\lambda}(t)$, $-\bar{\lambda}(t)$ vanishes. In
either case, $\lambda(t)=0$. Since this holds for all $\lambda$ such
that $u_{\lambda}\not=0$, and these span
$\mathfrak{c}^{*}\otimes\bc$, $t=0$.

\begin{Prop}
\label{convex}
Let $G$ be a semisimple Lie group. Let $\Gamma$ be a surface group. Let $\phi:\Gamma\to G$ be a homomorphism whose Zariski closure is reductive. Let $C$ denote the center of the centralizer of $\phi(\Gamma)$ in $G$. Let $P\subset\Lambda$ denote the subset of pure imaginary roots for which the quadratic form $\Im m(Q_{\lambda})$ is positive definite. Denote by $N=\Lambda\setminus\pm P$ the set of roots which are neither in $P$ nor opposites of elements of $P$. The following are equivalent.
\begin{enumerate}
  \item The central quadratic model has at least one smooth point.
  \item The central quadratic model has a dense set of smooth points.
  \item The convex hull of $\Im m(P))+\mathrm{span}(\Re e(N)\cup\Im m(N))$ in $\mathfrak{c}^*$ contains $0$ in its interior.
\end{enumerate}
\end{Prop}

\begin{pf}
Let $P'\subset \Lambda$ be a set which contains for each root
$\lambda\in N=\Lambda\setminus\pm P$ a representative of the pair
$\{\pm\lambda\}$ (resp. of the 4-uple
$\{\pm\lambda,\pm\bar{\lambda}\}$). Split $P'=R\cup I\cup M$
according to whether $\lambda$ is real, pure imaginary or neither.
When $\lambda\in R$, let $C_{\lambda}=t_{\lambda}$ and
$q_{\lambda}=Q_{\lambda}$. When $\ell\in P\cup I$, let
$C_{\lambda}=-\Im m(t_{\lambda})$ and $q_{\lambda}=Q_{\lambda}$. For
$E\subset \Lambda$, denote by $C_E=\{C_{\lambda}\,|\,\lambda\in
E\}$. With this notation, we get an expansion for $Q^{\mathfrak{c}}$
in real terms only,
\begin{eqnarray*}
Q^{\mathfrak{c}}=\sum_{\lambda\in P\cup I\cup R}q_{\lambda}C_{\lambda}+\sum_{\lambda\in M} \Re e(Q_{\lambda})\Re e(t_{\lambda})-\Im m(Q_{\lambda})\Im m(t_{\lambda})
\end{eqnarray*}

Assume that $\bar{\mathcal{Q}}^{\mathfrak{c}}_{\phi}$ contains at least one smooth point. Since this set is a cone, smooth points can be found arbitrarily close from $0$.  According to Lemma \ref{smoothcenter}, this means
that there exists a small nonzero $v=\sum_{\lambda\in \{0\}\cup
P\cup P'}v_{\lambda}$, $v_{\lambda}\in
H^1(\Gamma,\mathfrak{g}_{\lambda,\br})$, such that
$Q^{\mathfrak{c}}(v)=0$, and such that the set $E$ of roots
$\lambda\in\Lambda$ such that $v_{\lambda}\not=0$ spans
$\mathfrak{c}^{*}\otimes\bc$. {\em A fortiori}, $(E\cap P)\cup P'$
spans $\mathfrak{c}^{*}\otimes\bc$, i.e. $C_{E\cap P}\cup C_I \cup
C_R \cup \Re e(t_{M})\cup \Im m(t_M)$ spans $\mathfrak{c}$. Write
\begin{eqnarray*}
0=Q(v)=\sum_{\lambda\in E\cap P}q_{\lambda}(v_{\lambda})C_{\lambda}
+\sum_{\lambda\in I\cup R}q_{\lambda}(v_{\lambda})C_{\lambda}+\sum_{\lambda\in M}\Re e(Q_{\lambda}(v_{\lambda})t_{\lambda}).
\end{eqnarray*}
By construction, for $\lambda\in E\cap P$,
$q_{\lambda}(v_{\lambda})>0$. Thus 0 of $\mathfrak{c}$ belongs to
the interior of the convex hull of the image of $C_{E\cap P}$ in the
quotient vector space $\mathfrak{c}/\mathrm{span}(C_{R\cup I}\cup
\Re e(t_{M})\cup \Im m(t_M))$. In other words, 0 of $\mathfrak{c}$
belongs to the interior of the sum of the convex hull of $C_{E\cap
P}$ and the linear span of $C_{R\cup I} \cup \Re e(t_{M})\cup \Im
m(t_M)$. This implies that
$0\in\mathrm{interior}(\mathrm{convex\,hull}(\Im
m(P))+\mathrm{span}(\Re e(N)\cup\Im m(N)))$.

Conversely, assume that 0 of $\mathfrak{c}$ belongs to the interior
of the sum of the convex hull of $\Im m(P)$ and the linear span of
the real and imaginary parts of elements of $N$. In other words,
there exists a subset $E\subset P\cup P'$, small positive numbers
$a_{\lambda}$ for $\lambda\in E\cap P$, small nonzero real numbers
$a_{\lambda}$ for $\lambda\in E\cap (R\cup I)$ and small nonzero
complex numbers $c_{\lambda}$ for $\lambda\in E\cap M$ such that
\begin{eqnarray*}
\sum_{\lambda\in (E\cap P)\cup (E\cap (R\cup
I))}a_{\lambda}C_{\lambda}+\sum_{\lambda\in M}\Re
e(c_{\lambda}t_{\lambda})=0.
\end{eqnarray*}
For $\lambda\notin E$, let $a_{\lambda}=0$ or $c_{\lambda}=0$. Note
that  $span_{\br}(C_{(E\cap P)\cup (E\cap (R\cup I)))}\cup\Re
e(t_{E\cap M})\cup \Im m(t_{E\cap M}))=\mathfrak{c}$, thus
$span_{\bc}(E)=\mathfrak{c}^{*}\otimes\bc$.

Let $u\in H^{1}(\Gamma,\mathfrak{g})$ satisfy
$Q^{\mathfrak{c}}(u)=0$.  For $\lambda\in P\cup R\cup I$, there
exists a small $v_{\lambda}\in
H^{1}(\Gamma,\mathfrak{g}_{\lambda,\br})$ such that
$Q_{\lambda}(u_{\lambda}+v_{\lambda})=Q_{\lambda}(u_{\lambda})+a_{\lambda}$.
Indeed, for $\lambda\in E\cap P$, $q_{\lambda}$ is positive
definite, and takes all positive values. For $\lambda\in R\cup I$,
$q_{\lambda}$ is indefinite, and takes all real values. In the same
manner, for $\lambda\in M$, $Q_{\lambda}$ is $\bc$-bilinear and
nondegenerate, and takes all complex values, thus there exists a
small $c_{\lambda}$ such that
$Q_{\lambda}(u_{\lambda}+v_{\lambda})=Q_{\lambda}(u_{\lambda})+c_{\lambda}$.
For $\lambda\in P\setminus E$, take $v_{\lambda}=0$. Let
$v=\sum_{\lambda\in P\cup P'}v_{\lambda}$. Then
\begin{eqnarray*}
Q^{\mathfrak{c}}(u+v)
&=&\sum_{\lambda\in P\setminus E}q_{\lambda}(u_{\lambda})C_{\lambda}+\sum_{\lambda\in (E\cap P)\cup (E\cap (R\cup I))}q_{\lambda}(u_{\lambda}+v_{\lambda})C_{\lambda}\\
&&+\sum_{\lambda\in M}\Re e(Q_{\lambda}(u_{\lambda}+v_{\lambda})t_{\lambda})\\
&=&\sum_{\lambda\in P\setminus E}q_{\lambda}(u_{\lambda})C_{\lambda}+\sum_{\lambda\in (E\cap P)\cup (E\cap (R\cup I))}(q_{\lambda}(u_{\lambda})+a_{\lambda})C_{\lambda}\\
&&+\sum_{\lambda\in M}\Re e(Q_{\lambda}(u_{\lambda})+c_{\lambda})t_{\lambda})\\
&=&Q^{\mathfrak{c}}(u)+\sum_{\lambda\in (E\cap P)\cup (E\cap (R\cup I))}a_{\lambda}C_{\lambda}+\sum_{\lambda\in M}\Re e(c_{\lambda}t_{\lambda})\\
&=&0.
\end{eqnarray*}
By construction, the set of roots such that $(u+v)_{\lambda}\not=0$
contains $E$,  it spans $\mathfrak{c}^{*}$, thus, with Lemma
\ref{smoothcenter}, $u+v$ gives a smooth point of the central
quadratic model (in cohomology). Hence smooth points are dense in the central
quadratic model.
\end{pf}

\subsection{Flexibility for toral centralizers}
\label{flextor}

Here comes the proof of Theorem \ref{flexssimple} in case the centralizer is abelian.

\begin{co}
\label{abelian} Let $G$ be a semisimple real algebraic group. Let
$\Gamma$ be the fundamental  group of a surface of genus $\geq
\mathrm{dim}(G)^{2}$. Let $\phi:\Gamma\to G$ be a reductive
homomorphism. Assume that the centralizer $\mathfrak{c}$ of
$\phi(\Gamma)$ is abelian. Split $\mathfrak{g}$ into real root
spaces $\mathfrak{g}_{\lambda,\br}$ under $\mathfrak{c}$. When
$\lambda$ is pure imaginary, this is a symplectic vector space,
which gives rise to a Toledo invariant $T_{\lambda}$. For all
nonzero pure imaginary roots $\lambda$, the following Milnor-Wood
type inequality holds.
\begin{eqnarray}
\label{MW}
4|T_{\lambda}|\leq -\chi(\Gamma)\mathrm{dim}(\mathfrak{g}_{\lambda,\br}).
\end{eqnarray}
Denote by $P$ denotes the set of nonzero pure imaginary roots $\lambda\in\Lambda$ such that $4T_{\lambda}= -\chi(\Gamma)\mathrm{dim}(\mathfrak{g}_{\lambda,\br})$. Denote by $N$ the set of roots which neither in $P$ nor opposites of elements of $P$. Then $\phi$ is flexible if and only if
\begin{eqnarray*}
0\in \mathrm{interior}(\mathrm{convex\,hull}(\Im m(P))+\mathrm{span}(\Re e(N)\cup\Im m(N))).
\end{eqnarray*}
\end{co}

\begin{pf}Note that
$H^0(\Gamma,\mathfrak{g}_{\lambda,\br})=H^2(\Gamma,\mathfrak{g}_{\lambda,\br})=0$,
hence $\mathrm{dim}
H^1(\Gamma,\mathfrak{g}_{\lambda,\br})=-\chi(\Gamma)\mathrm{dim}(\mathfrak{g}_{\lambda,\br}).$
The signature of a quadratic form on a vector space is less than or
equal to the dimension of this space, this gives inequality
(\ref{MW}). Equality holds if and only if $Q_{\lambda}$ is definite.
Thus the set $P$ of roots coincides with the set defined in
Proposition \ref{convex}. Proposition \ref{convex} yields the
criterion for the quadratic model $\mathcal{Q}_{\phi}$ to be smooth
of dimension $-\chi(\Gamma)\mathrm{dim}(G)$ at a dense set of
points. Theorem \ref{dgms} implies that the same holds for $Hom(\Gamma,G)$
in a neighborhood of the conjugacy class of $\phi$.
\end{pf}

\begin{ex}
\label{su(1,1)->su(2,1)} Let $\Gamma$ be a cocompact lattice in
$H=SU(1,1)$. Embed $SU(1,1)$ as a lower right diagonal block in
$G=SU(2,1)$. Then the obtained homomorphism $\Gamma\to SU(2,1)$ is
not flexible.
\end{ex}

In this case, treated originally by W. Goldman, \cite{G1},  the
centralizer $\mathfrak{z}$ of $H$ is the 1-dimensional subspace
generated by $Z=diag(-2i,i,i)$. Under $ad_Z$, $\mathfrak{g}=
\{\left[\begin{matrix}
       v_1 & v_2  &  v_3 \\
     -\bar{v}_2 & v_4  &  v_5 \\
      \bar{v}_3& \bar{v}_5 & v_6 \end{matrix}\right ]\}$, with $v_2$, $v_3$, $v_5 \in \bc$ and
       $v_1$, $v_4$, $v_6 \in\Im m(\bc)$,  splits as
\begin{eqnarray*}
 \mathfrak{g}=\mathfrak{g}_0 \oplus \mathfrak{g}_{3i,\br},
\end{eqnarray*}
where   $\mathfrak{g}_0 =\mathfrak{z}\oplus \mathfrak{su}(1,1)$ and
$\mathfrak{g}_{3i,\br}=\{\left[\begin{matrix}
       0 & v_2  &  v_3 \\
     -\bar{v}_2 & 0  &  0 \\
      \bar{v}_3& 0 & 0 \end{matrix}\right ]\}$ is a 4-dimensional real vector space
isomorphic to $\bc^2$ with $Z$ acting by multiplication by
$diag(3i,3i)$ and $SU(1,1)$ by its standard representation.
The restriction of the Killing form to $\mathfrak{g}_{3i,\br}$ is a
quadratic form $q$ of signature $(1,1)$ since $v_2\in \mathfrak{t}$
and $v_3\in \mathfrak{p}$ for the Cartan decomposition
$\mathfrak{g}=\mathfrak{t}\oplus\mathfrak{p}$.

Consider the flat bundle $E_{3i}$ on $\Sigma$ with fiber
$\mathfrak{g}_{3i,\br}$ equipped with  the symplectic form $\Im m
\Omega_{3i}$. View the lifted bundle
$\tilde{E_{3i}}=\tilde{\Sigma}\times\mathfrak{g}_{3i,\br}$ as a
trivial bundle over $H^{1}_{\bc}=\tilde{\Sigma}$. Let $\tilde{\tau}\subset
\tilde{E_{3i}}$ denote the tautological complex line bundle, i.e.,
at a point $p\in H^{1}_{\bc}$ representing a negative line $\ell$, the fiber
of $\tilde{\tau}$ is $\ell$. Let $\tilde{\tau}^{\bot}$ be its
orthogonal complement. Both descend to bundles over $\Sigma$, and
$E_{3i}=\tau\oplus\tau^{\bot}$. The quadratic form $q$ is negative
definite on $\tau$ and positive definite on $\tau^\bot$. Let $J$ be
the (variable) complex structure on $W$ defined by $J=i$ on $\tau$
and $J=-i$ on $\tau^\bot$. Then for $B=\begin{pmatrix} b_-\\ b_+
\end{pmatrix}\in E_{3i}=\tau\oplus\tau^{\bot}$, $\Im m\Omega_{3i}(B,JB)=-q(b_-)+q(b_+)$ is positive definite.
Thus $J$ is tamed by $\Im m\Omega_{3i}$. As a complex plane bundle,
$(E_{3i},J)$ is isomorphic to $E_{3i}^{1,0}=\{B\in
E_{3i}\otimes\bc\,|\,JB=iB\}=\tau\otimes\bc\simeq \tau\oplus\tau$.
Therefore $c_1 (E_{3i},\Omega_{3i})=c_1
(\tau\oplus\tau)=2c_1(\tau)$. Recall that, as a homogeneous vector
bundle, the tangent bundle of $H^{1}_{\bc}$ is
$TH^{1}_{\bc}=Hom(\tilde{\tau},\tilde{\tau}^\bot)=\tilde{\tau}^*\otimes\tilde{\tau}^\bot$.
Therefore, as bundles on $\Sigma$, $T\Sigma=\tau^*\otimes\tau^\bot$,
$c_1 (T\Sigma)=c_1 (\tau^*)+c_1(\tau^\bot)$. Since
$\tau\oplus\tau^\bot =E_{3i}$ is flat, $c_1 (\tau^\bot)=-c_1(\tau)$,
and $c_1 (T\Sigma)=-2c_1 (\tau)$. Therefore $c_1
(E_{3i},\Omega_{3i})=-c_1 (T\Sigma)$. As a number, $c_1
(E_{3i},\Omega_{3i})=-\chi(\Gamma)$. So the Milnor-Wood inequality
is an equality.

\subsection{Closedness of flexibility}

In principle, Corollary \ref{abelian} gives a complete description of a neighborhood of a homomorphism in $Hom(\Gamma,G)$, when its centralizer is abelian. It follows that

\begin{co}
\label{openclosed}
Let $G$ be a semisimple Lie group. Let $\Gamma$ be the fundamental group of a closed surface of genus $g\geq \mathrm{dim}(G)^{2}$. Let $T$ denote a (non necessarily maximal) torus in $G$. Let $\chi_{T}$ denote the set of conjugacy classes of reductive homomorphisms $\Gamma\ra G$ whose centralizers are conjugate to $T$. The set of flexible conjugacy classes in $\chi_{T}$ is both open and closed in $\chi_{T}$.
\end{co}

\begin{pf}
Up to a linear change of coordinates, the central quadratic model
$\mathcal{Q}_{\phi}^{\mathfrak{c}}$ is determined by $T$ and by a
collection of Toledo invariants. This is because once $T$ is fixed,
$\mathfrak{t}_\lambda$ is determined, and $Q_\lambda$ is determined
by its signature, hence by the first Chern class up to the linear
change of coordinates. Here,
$\mathfrak{z}=\mathfrak{t}=\mathfrak{c}$. Since Toledo invariants
are locally constant on $Hom(\Gamma,G)$, so is the quadratic model
$\mathcal{Q}_{\phi}/T$ as $\phi$ moves in $\chi_{T}$.

An alternate proof consists in invoking Theorem \ref{thm1}.
\end{pf}

\begin{ex}
\label{sl2block}
Let $\phi_1$, $\phi_2 :\Gamma\ra Sl(2,\br)$ be two nonconjugate discrete cocompact homomorphisms. Use them to map $\Gamma$ as a Zariski-dense subgroup of $Sl(2,\br)\times Sl(2,\br)$. Then map $Sl(2,\br)\times Sl(2,\br)$ to $Sl(4,\br)$. The obtained homomorphism is flexible.
\end{ex}
Indeed, the centralizer of $Sl(2,\br)\times Sl(2,\br)$ in $Sl(4,\br)$ is 1-dimensional, so the set $\mathcal{P}$ of considered homomorphisms consists in a connected component of a single $\chi_{T}$. $\mathcal{P}$ contains in its closure the diagonal embedding studied in Example \ref{exsscentr}, which is flexible. Since flexibility is an open condition, $\mathcal{P}$ contains flexible elements, thus all of $\mathcal{P}$ is flexible, according to Corollary \ref{openclosed}.

\section{Levi factors of centralizers}
\label{levi}

\subsection{The trivial homomorphism}
\label{trivialh}
The simplest case of a homomorphism whose centralizer has a nontrivial Levi factor is the trivial homomorphism. This special case is of course classical. The method and some intermediate results will turn out to be useful later on.

\begin{lemma}
\label{trivial}
Let $\Sigma$ be a closed surface with genus $g\geq 1$. Let $\Gamma=\pi_1 (\Sigma)$ act trivially on a semisimple Lie algebra $\mathfrak{g}$. Let $Q$ denote the quadratic map $H^1 (\Gamma,\mathfrak{g})\to H^2 (\Gamma,\mathfrak{g})$, $u\mapsto [u\smile u]$. Let us view an element $u\in H^1 (\Gamma,\mathfrak{g})$ as a linear map $H_{1}(\Gamma,\br)\ra\mathfrak{g}$. Define $span(u)$ as the image of this map. If the centralizer of $span(u)$ in $\mathfrak{g}$ is trivial, then
\begin{enumerate}
\item $Q$ is a submersion at $u$ ;
\item the $G$-action on $H^1 (\Gamma,\mathfrak{g})$ is locally free at $u$.
\end{enumerate}
\end{lemma}

\begin{pf}
1. Denote by $(X,Y)\mapsto X\cdot Y$ the Killing form on $\mathfrak{g}$. Assume $Q:u\mapsto [u\smile u]$ is not a submersion at $u_0$. This means that the linear map $v\mapsto [u_0 \smile v]$, $H^1 (\Gamma,\mathfrak{g})\ra H^2 (\Gamma,\mathfrak{g})=\mathfrak{g}$, is not onto. Let $X\in\mathfrak{g}$ be orthogonal to its image. Let $Y_1, \ldots,Y_k$ be a basis of $span(u_0)$. Then
\begin{eqnarray*}
u_0 =\sum_{i=1}^{k}a_i \otimes Y_i ,
\end{eqnarray*}
where $a_i \in H^1 (\Gamma,\br)$ are linearly independent. By Poincar\'e duality, the pairing $(a,b)\mapsto u\smile v$ on $H^1 (\Gamma,\br)$ is nondegenerate. Therefore, there exist classes $b_i \in H^1 (\Gamma,\br)$ such that $a_i \smile b_j =\delta_{ij}$. For $Y\in\mathfrak{g}$, take $v_{i,Y}=b_i \otimes Y$. Then $[u_0 \smile v_{i,Y}]=[Y_i ,Y]$,
$$0=X\cdot[u_0 \smile v_{i,Y}]=X\cdot[Y_i ,Y]=-[Y_i ,X]\cdot Y.$$
This shows that for all $i=1,\ldots,k$, $[Y_i ,X]=0$, i.e. $X$ centralizes the subalgebra generated by $span(u_0)$.

2. The differential of the $G$-action on $H^1 (\Gamma,\mathfrak{g})$ is $(Z,u)\mapsto ad_{Z}(u)$, for $Z\in \mathfrak{g}$. If $u_{0}$ is viewed as a linear map $H_{1}(\Gamma,\br)\ra\mathfrak{g}$, then $ad_{Z}(u_{0})=ad_{Z}\circ u_{0}$. This vanishes only if $Z$ centralizes the image $span(u_{0})$ of $u_{0}$.
\end{pf}

\begin{co}
\label{cortrivial} Let $G$ be a semisimple Lie group. Let $\Gamma$ be the fundamental group of a closed surface of genus $>1$. The trivial homomorphism of $\Gamma$ to $G$ belongs to the closure of $\mathcal{S}$.
\end{co}

\begin{pf}
Take two vectors $X$ and $Y$ in $\mathfrak{g}$ which do not have a common centralizer (this is an open dense condition). Take non colinear classes $a$ and $b$ in $H^1 (\Gamma,\br)$ such that $a\smile b=0$. Then $v =a\otimes X+b\otimes Y$ satisfies $[v \smile v]=2(a\smile b)[X,Y]=0$. Since $v$ does not belong to $H^1 (\Gamma,\br)\otimes\mathfrak{z}$ for any centralizer $\mathfrak{z}$, Lemma \ref{trivial} ensures that the quadratic cone $\{u\,|\,[u\smile u]=0\}$ is smooth at $v$. Furthermore, since $X$ and $Y$ do not have a common centralizer, the $G$-action is locally free near $v$. One concludes that the quadratic model $\bar{\mathcal{Q}}_{triv}=\{u\,|\,[u\smile u]=0\}$ is smooth near the equivalence class of $v$, which can be chosen arbitrarily close to the origin.
\end{pf}

\begin{rem}
This proves Theorem \ref{thm0} for semisimple Lie groups. The generalization to reductive groups is postponed until subsection \ref{proofthm1}.
\end{rem}

Note that the preceding argument is not sufficient to prove flexibility of the trivial representation itself, i.e. density of smooth homomorphisms in a neighborhood of the trivial homomorphism. Instead of constructing just one explicit class $v$, one needs to adapt $v$ to a given solution $u$ of $[u\smile u]=0$. We start with some linear algebraic preliminaries.

\begin{lemma}
\label{grass} Let $(V,\Omega)$ be a $2g$-dimensional symplectic
vector space. Let $W$ be a $n$-dimensional vector space, and
$\omega$ an alternating 2-form on $W$. Let $E$ be a codimension $q$
vectorsubspace in $Hom(W,V)$. Let $f_{0}$ be an element of
$Hom(W,V)$. If $g\geq 2nq$, there exists a linear map $f:W\ra V$
such that
\begin{enumerate}
\item $f\in E$ ;
\item $f_{0}+\epsilon f$ is injective for all nonzero $\epsilon$ ;
\item $f^{*}\Omega=\omega$ ;
\item the map $\Phi_{f}:E\to W^{*}\otimes W^{*}$ defined for $w$, $w'\in W$ by
$$\Phi_{f}(f')(w,w')=\Omega(f(w),f'(w'))$$
is onto.
\end{enumerate}
\end{lemma}

\begin{pf}
First, we prove that there exists a symplectic $2n$-dimensional
subvector space $P\subset V$ such that $Hom(W,P)\subset E$.

Let $\mathcal{G}$ be the Grassmannian of $2n$-planes in $V$, and let $\mathcal{S}$ denote the subset of symplectic planes, i.e. those on which $\Omega$ induces a symplectic structure. Let $\tau$ denote the tautological $2n$-bundle on $\mathcal{G}$.

Let $\ell$ be a linear form on $Hom(W,V)$. Then, for each $P\in \mathcal{G}$, $Hom(W,P)$ is a subspace of $Hom(W,V)$, therefore $\ell$ defines a linear form on $Hom(W,P)$, i.e. an element of $Hom(P,W)$. In other words, $\ell$ defines a section of the bundle $Hom(\tau,W)$, which is isomorphic to the direct sum of $n$ copies of $\tau^{*}$. Since $\tau^{*}$ is isomorphic to $\tau$, $\ell$ defines as well a section of $\tau^{n}$, the direct sum of $n$ copies of $\tau$. Choosing $q$ independant linear equations for $E$ defines a section $s$ of $\tau^{nq}$.

Let us show that $s$ must vanish somewhere on $\mathcal{S}$.  For
this, since $\tau^{nq}_{|\mathcal{S}}$ is oriented, it is sufficient
to prove that the Euler class $e(\tau^{nq}_{|\mathcal{S}})$ does not
vanish. Pick a complex structure $J$ on $V$ compatible with $\Omega$
and let $\mathcal{C}$ denote the Grassmannian of complex
$n$-dimensional $J$-complex subspaces in $V$. Since
$\mathcal{C}\subset\mathcal{S}$, it suffices to show that
$e(\tau^{nq}_{|\mathcal{C}})$ does not vanish. As a real
vectorbundle, $\xi=\tau_{|\mathcal{C}}$ coincides with the
tautological complex $n$-plane bundle $\tau_{\bc}$ over
$\mathcal{C}$. In other words,
$\xi\otimes\bc=\tau_{\bc}\oplus\overline{\tau_{\bc}}$. Thus
\begin{eqnarray*}
e(\xi)^2 &=&e(\xi\oplus\xi)\\
&=&(-1)^{n}e(\xi\otimes\bc)\\
&=&(-1)^{n}c_{n}(\tau_{\bc})c_{n}(\overline{\tau_{\bc}})\\
&=&c_{n}(\tau_{\bc})^2 ,
\end{eqnarray*}
and
\begin{eqnarray*}
e(\tau^{2nq}_{|\mathcal{C}})=e(\xi)^{2nq}=c_{n}(\tau_{\bc})^{2nq}.
\end{eqnarray*}
Let $\mathcal{C}_{\infty}$ denote the direct limit of the
Grassmannians $\mathcal{C}$ as $g$ tends to infinity.  According to
A. Borel, \cite{Borel_coh}, the cohomology algebra of
$\mathcal{C}_{\infty}$ is freely generated by the Chern classes of
the tautological bundle $\tau_{\bc,\infty}$. Thus
$c_{n}(\tau_{\bc,\infty})^{2nq}\not=0$. It follows that
$e(\tau^{nq}_{|\mathcal{C}})$ does not vanish for $g$ large enough.
In fact, $g\geq 2nq$ is sufficient, see \cite{Manivel}, chapter 3,
section 2. We conclude that for $g\geq 2nq$, there exists a
$J$-complex $n$-plane $P\in V$ such that for all $\br$-linear maps
$f:W\to P$, $f\in E$.

Let $r$ denote the rank of $\omega$. Since $r\leq n$,  there exists
an $n$-dimensional subspace $P'\subset P$, transverse to
$\mathrm{Im}(f_{0})\cap P$, such that $\Omega_{|P'}$ has rank $r$.
Let $f_1 : W\ra P'$ be a bijective linear map. Then the pulled-back
form $f_{1}^{*}\Omega$ has rank $r$. There exists a bijective linear
map $h:W\ra W$ such that $h^{*}(f_{1}^{*}\Omega)=\omega$. Then
$f=f_1 \circ h$ satisfies $f^{*}\Omega=\omega$ and $f$ belongs to
$E$. $f$ is injective, $\mathrm{Im}(f)$ and $\mathrm{Im}(f_{0})$ are
transverse, thus $f_{0}+\epsilon f$ is injective for all nonzero
$\epsilon$.

Let $\phi$ denote the restriction of $\Phi_f$ to $Hom(W,P)\subset
E$. Then
\begin{eqnarray*}
\mathrm{ker}(\phi)&=&\{f'\in Hom(W,P)\,|\,\forall w,\,w'\in W,\,\Omega(f(w),f'(w'))=0\}\\
&=&\{f'\in Hom(W,P)\,|\,\forall w,\,w'\in W,\,f'(w')\bot P'\}\\
&=& Hom(W,P'^{\bot}).
\end{eqnarray*}
In particular, $\mathrm{dim}(\mathrm{ker}(\phi))=n^2$,
$\mathrm{dim}(\mathrm{im}(\phi))=2n^2 -n^2 =n^2$,  i.e. $\phi$ (and
thus $\Phi_f$) is onto.
\end{pf}

\begin{lemma}
\label{trivialbis} Let $\Sigma$ be a closed surface with genus
$g\geq 1$. Let $\Gamma=\pi_1 (\Sigma)$ act trivially on  a
semisimple Lie algebra $\mathfrak{g}$. Let $Q$ denote the quadratic
map $H^1 (\Gamma,\mathfrak{g})\to H^2 (\Gamma,\mathfrak{g})$,
$u\mapsto [u\smile u]$. Let $L:H^1 (\Gamma,\mathfrak{g})\to H^2
(\Gamma,\mathfrak{g})$ be an arbitrary linear map. Let $v_{0}$ be an
arbitrary element of $H^1 (\Gamma,\mathfrak{g})$. Let
$\mathcal{U}_{L,v_{0}}\subset H^1 (\Gamma,\mathfrak{g})$ denote the
subset of elements $v$ such that
\begin{enumerate}
\item $L(v)=0$ ;
\item $v_{0}+\epsilon v$ has full span for all nonzero $\epsilon$ ;
\item the map
\begin{eqnarray*}
H^{1}(\Gamma,\mathfrak{g})\to H^2 (\Gamma,\mathfrak{g}),\quad u'\mapsto L(u')+[v\smile u']
\end{eqnarray*}
is onto.
\end{enumerate}
If $g\geq 2\mathrm{dim}(\mathfrak{g})^{2}$,  then
$Q:\mathcal{U}_{L,v_0}\to H^2 (\Gamma,\mathfrak{g})$ is onto.
\end{lemma}

\begin{pf}
Let $e_1,\ldots,e_{2g}$ be a basis of $H^1 (\Gamma,\br)$,
$Z_1,\cdots,Z_n$  a basis of $\mathfrak{g}$. If
$u=\sum_{i=1}^{2g}\sum_{j=1}^{n}x_{ij}e_i \otimes Z_j$,
$u'=\sum_{i=1}^{2g}\sum_{j=1}^{n}x'_{ij}e_i \otimes Z_j$, then
$[u\smile u']=\sum_{i,j,k,l}x_{ik}e_i\cup
e_jx'_{jl}[Z_k,Z_l]=\sum_{k,\,\ell=1}^{n}y_{k\ell}[Z_k ,Z_\ell]$,
where the matrices $X=(x_{ij})_{i=1,\ldots,2g,\,j=1,\ldots n}$,
$Y=(y_{k\ell})_{k,\,\ell=1,\ldots n}$ and $A=(e_i \smile
e_j)_{i,\,j=1,\ldots,2g}$ are related by
\begin{eqnarray*}
Y=X^{\top}AX'.
\end{eqnarray*}
In other words, the cup-product on $H^1 (\Gamma,\br)$ induces a bilinear map $(u,u')\ra u\smile u'$, $H^1 (\Gamma,\br)\otimes\mathfrak{g}\ra H^2 (\Gamma,\br)\otimes\mathfrak{g}\otimes\mathfrak{g}$, whose image by the Lie bracket $[\cdot,\cdot]:\mathfrak{g}\otimes\mathfrak{g}\to\mathfrak{g}$ equals $[u\smile u']\in H^2 (\Gamma,\mathfrak{g})$. Note that $u\smile u\in H^2 (\Gamma,\br)\otimes\Lambda^{2}\mathfrak{g}$.

Let $V=H^1 (\Gamma,\br)$ be equipped with the symplectic  form
$\Omega$ defined by cup-product evaluated on the fundamental class,
and $W=\mathfrak{g}^{*}$. If we identify $ H^1
(\Gamma,\mathfrak{g})$ with
$Hom(H_1(\Gamma,\br),\mathfrak{g})=Hom(V^*,\mathfrak{g})=\mathfrak{g}\otimes
V=Hom(W,V)$, one can identify $u\in H^1 (\Gamma,\mathfrak{g})$ with
a map $f_u\in Hom(W,V)$. Note the correspondence is given, for $w\in
\mathfrak{g}^*, \gamma\in \Gamma$, by
$$f_u(w)(\gamma)=w(u(\gamma)).$$

 Then, since the cup product in $H^1(\Gamma,\mathfrak{g})$ corresponds to the symplectic form $\Omega$ via above identification,
 using the notation of Lemma
\ref{grass}, $\Phi_{f_{u}}(f_{u'})(w,w')=\Omega(f_uw,f_{u'}w')$, and
hence $u\smile u$ identifies with $f_{u}^{*}\Omega$ and $u\smile u'$
with $\Phi_{f_{u}}(f_{u'})$. We define a subspace $E$ of
$Hom(W,V)=H^1 (\Gamma,\mathfrak{g})$ by
$$
f_u \in E\Leftrightarrow L(u)=0.
$$
As an element $f_{0}$ of $Hom(W,V)$, pick $f_{v_{0}}$.

Since $\mathfrak{g}$ is semi-simple, every $Z\in\mathfrak{g}=H^2 (\Gamma,\mathfrak{g})$ is the image under the Lie bracket of some $\omega\in\Lambda^{2}\mathfrak{g}$. Lemma \ref{grass} provides us with an element $v\in H^1 (\Gamma,\mathfrak{g})$ such that $L(v)=0$ and $v\smile v=\omega$, whence $Q(v)=Z$. Lemma \ref{grass} asserts that $\Phi_{f_{v}}:E\to W^* \otimes W^*$, is onto. This translates into surjectivity of
$$
\mathrm{ker}(L)\ra H^2 (\Gamma,\br)\otimes\mathfrak{g}\otimes\mathfrak{g}, \quad u'\mapsto v\smile u'.
$$
Again, this implies that the map $u'\mapsto [v\smile u']$, $\mathrm{ker}(L)\ra H^2 (\Gamma,\mathfrak{g})$ is onto. A fortiori,
\begin{eqnarray*}
H^{1}(\Gamma,\mathfrak{g})\ra H^2 (\Gamma,\mathfrak{g}),\quad u'\mapsto L(u')+[v\smile u']
\end{eqnarray*}
is onto. Note that $f_{v_{0}}+\epsilon f_v$ being injective implies
that $v_{0}+\epsilon v$ has full span. Indeed $((f_{v_{0}}+\epsilon
f_v)(w))(\gamma)=w((v_0+\epsilon v)(\gamma))$, hence if
$(f_{v_{0}}+\epsilon f_v)(w)=0$, then $w((v_0+\epsilon
v)(\gamma))=0$ for all $\gamma\in\Gamma$, so $v_{0}+\epsilon v$ must
have full span to have $w=0$.

Therefore $v\in\mathcal{U}_{L,v_{0}}$. We conclude that the image
$Q(\mathcal{U}_{L,v_{0}})$ contains all of $H^2
(\Gamma,\mathfrak{g})$.
\end{pf}

Lemma \ref{trivialbis} leads to the flexibility of the trivial representation. It is also the key to the proof of  Theorem \ref{flexssimple} in the next subsection.

\begin{co}
\label{flextriv}
Let $G$ be a semi-simple real algebraic group. Let $\Sigma$ be a compact surface with genus $\geq\mathrm{dim}(G)^{2}$. Then the trivial homomorphism $\Gamma\ra G$ is flexible.
\end{co}

\begin{pf}
Let $u_0 \in H^1 (\Gamma,\mathfrak{g})$ be such that $[u_0 \smile u_0]=0$. We need a $v\in H^1 (\Gamma,\mathfrak{g})$ such that $u_{\epsilon}=u_0 +\epsilon v$ has full span and satisfies $[u_{\epsilon}\smile u_{\epsilon}]=0$. We achieve this by requiring that $[u_0 \smile v]=0$ and $[v\smile v]=0$. Therefore we set $L(v)=[u_0 \smile v]$ and $v_0 =u_0$. Lemma \ref{trivialbis} yields a cohomology class $v$ such that $[u_0 \smile v]=0$, $[v\smile v]=0$ and $u_0 +\epsilon v$ has full span for all nonzero $\epsilon$. According to Lemma \ref{trivial}, the conjugacy class of $u_0 +\epsilon v$ is a smooth point of the quadratic model $\bar{\mathcal{Q}}_{triv}$. With Theorem \ref{dgms} and Corollary \ref{dense}, this shows that Zariski dense representations are dense in a neighborhood of the trivial representation.
\end{pf}

\subsection{Semisimple centralizers}
\label{proofrigidssimple}

Here comes the proof of an other special case of Theorem
\ref{flexssimple}: when the centralizer is semisimple.

\begin{Prop}
\label{surface1}
Let $G$ be a semi-simple real algebraic group. Let $\Sigma$ be a compact surface with genus $\geq 2\mathrm{dim}(G)^{2}$ and $\phi:\pi_1(\Sigma)=\Gamma \ra G$ be a homomorphism with reductive Zariski closure. Assume that its centralizer $\mathfrak{z}$ is semi-simple. Then $\phi$ is flexible.
\end{Prop}

\begin{pf}
Let $H$ denote the Zariski closure of $\phi(\Gamma)$ and $Z=Z_{G}(H)$ its centralizer. Since $HZ$ is reductive, there exists an $HZ$-invariant splitting
$$
\mathfrak{g}=\mathfrak{z}\oplus\mathfrak{z}'.
$$

Let $u_0 \in H^{1}(\Gamma,\mathfrak{g})$ be such that $[u_0 \smile
u_0]=0$. Write $u_{0}=v_{0}+w_{0}$  with $v_{0}\in
H^1(\Gamma,\mathfrak{z})$ and $w_{0}\in H^1(\Gamma,\mathfrak{z}')$.
For $v\in H^1(\Gamma,\mathfrak{z})$, let $L(v)=[u_0 \smile v]\in H^2
(\Gamma,\mathfrak{z})=H^2 (\Gamma,\mathfrak{g})$.

Note that $\Gamma$ acts trivially on $\mathfrak{z}$. According to
Lemma \ref{trivialbis},  there exists $v\in
H^1(\Gamma,\mathfrak{z})$ such that $L(v)=0$ and $[v\smile v]=0$.
Then $u=u_0 +v$ satisfies $Q(u)=0$. The map
$$
u'\mapsto[u\smile u']=L(u')+[v\smile u'],\quad H^1(\Gamma,\mathfrak{z})\to H^2 (\Gamma,\mathfrak{z})=H^2 (\Gamma,\mathfrak{g})
$$
is onto, thus the cup-product map on $H^1(\Gamma,\mathfrak{g})$ is a submersion at $u$.

The $\mathfrak{z}$-component of $u$ is $v_{0}+v$. Since $v_{0}+v$ has full span, its stabilizer under $\mathfrak{z}$ is trivial. It follows that the $Z=Z_{G}(\phi(\Gamma))$-action on $H^1(\Gamma,\mathfrak{g})$ is locally free at $u$.

Since $v$ can be chosen to be arbitrarily small, one concludes that there are smooth points of the quadratic model $\bar{\mathcal{Q}}_{\phi}=\{[u \smile u]=0\}$ in all neighborhoods of the class of $u$. With Theorem \ref{dgms}, this shows that smooth points are dense in a neighborhood of $\phi$ in $Hom(\Gamma,G)$.
\end{pf}

\begin{ex}
\label{exsscentr}
Map $\Gamma$ to a discrete cocompact subgroup of $Sl(2,\br)$, then map $Sl(2,\br)$ diagonally to $Sl(2,\br)\times Sl(2,\br)$, and finally realize $Sl(2,\br)\times Sl(2,\br)$ as block-diagonal matrices in $Sl(4,\br)$. The obtained homomorphism is flexible.
\end{ex}
Indeed, the centralizer of the diagonal $Sl(2,\br)$ in $Sl(4,\br)$ is isomorphic to $Sl(2,\br)$.

\subsection{Reductive centralizers}
\label{reductivecent}

We show that occurrence of a nontrivial Levi factor in a centralizer can only bring extra flexibility, without distroying the flexibility allowed by the center of the centralizer.

\begin{Prop}
\label{extra}
Let $G$ be a semisimple real algebraic group. Let $\Gamma$ be the fundamental group of a closed surface $\Sigma$ of genus $\geq 2\mathrm{dim}(G)^{2}$. Let $\phi:\Gamma\to G$ be a homomorphism with reductive Zariski closure. Assume that the central local quadratic model $\mathcal{Q}_{\phi}^{\mathfrak{c}}$ for $Hom(\Gamma,G)$ near $\phi$ has at least one smooth point. Then the full quadratic model $\mathcal{Q}_{\phi}$ has a dense set of smooth points.
\end{Prop}

\begin{pf}
Let $H$ denote the Zariski closure of $\phi(\Gamma)$ and $Z=Z_{G}(H)$ its centralizer. Since $HZ$ is reductive, there exists an $HZ$-invariant splitting
$$
\mathfrak{g}=\mathfrak{z}\oplus\mathfrak{z}'.
$$

The Lie algebra of $Z$ splits as $\mathfrak{z}=\mathfrak{s}\oplus\mathfrak{c}$ where $\mathfrak{s}$ is semi-simple and $\mathfrak{c}$ is the center of $\mathfrak{z}$. Let $u\in H^{1}(\Gamma,\mathfrak{g})$ be such that $Q(u)=0$.

Assume that $\bar{\mathcal{Q}}_{\phi}^\mathfrak{c}$ has at least one
smooth point.  There exist arbitrarily small $v\in
H^{1}(\Gamma,\mathfrak{g})$ such that $[(u+v)\smile
(u+v)]^{\mathfrak{c}}=0$. Furthermore, one can assume that
\begin{itemize}
  \item the linearization of $Q^{\mathfrak{c}}$ at $u+v$, $u'\mapsto [(u+v)\smile u']^{\mathfrak{c}}$, is onto ;
  \item for $Z\in \mathfrak{c}$, $ad_{Z}(u+v)=0$ implies $Z=0$.
\end{itemize}

The arguments in the proof of Proposition \ref{surface1} take care
of the $\mathfrak{s}$-component $Q^\mathfrak{s}$ by taking a linear
map $L:H^1(\Gamma,\mathfrak{s})\ra H^2(\Gamma,\mathfrak{s})$,
$w\mapsto [(u+v)^{\mathfrak{s}}\smile w]$ and $(u+v)^{\mathfrak{s}}$
as an arbitrarily given element. Lemma \ref{trivialbis} provides us
with a small $w\in H^{1}(\Gamma,\mathfrak{s})$ such that
$[(u+v)^{\mathfrak{s}}\smile w]=0$ and $[w\smile
w]=-[(u+v)\smile(u+v)]^{\mathfrak{s}}$. Furthermore, one can assume
that $u^{\mathfrak{s}}+v^{\mathfrak{s}}+w$ has full span in
$\mathfrak{s}$, and that the map
\begin{eqnarray*}
H^{1}(\Gamma,\mathfrak{s})\ra H^{2}(\Gamma,\mathfrak{s}),\quad u'\mapsto [(u^{\mathfrak{s}}+v^{\mathfrak{s}})\smile u']+[w\smile u']
\end{eqnarray*}
is onto. Then
\begin{eqnarray*}
[(u+v+w)\smile (u+v+w)]&=&[(u+v+w)\smile (u+v+w)]^{\mathfrak{c}}\\
&&+[(u+v+w)\smile
(u+v+w)]^{\mathfrak{s}}
\end{eqnarray*}
\begin{eqnarray*}
&=&2[(u+v)\smile w]^{\mathfrak{c}}+[(u+v)\smile (u+v)]^\mathfrak{s}+[w\smile w]\\
&=&2[((u+v)^\mathfrak{z}+(u+v)^\mathfrak{z'})\smile w]^\mathfrak{c}
=2[(u+v)^\mathfrak{z}\smile w]^\mathfrak{c}\\
&=& 2[((u+v)^\mathfrak{c}+(u+v)^\mathfrak{s})\smile w]^\mathfrak{c}=2[    (u+v)^\mathfrak{s}\smile w]^\mathfrak{c}=0.
\end{eqnarray*}
Also $u'\in H^1(\Gamma,\mathfrak{s})\mapsto [(u+v+w)\smile
u']^\mathfrak{s}=[(u+v)^\mathfrak{s}\smile u']+[w\smile u']$ is onto
to $H^2(\Gamma,\mathfrak{s})$ and $u'\in
H^1(\Gamma,\mathfrak{c})\mapsto [(u+v+w)\smile
u']^\mathfrak{c}=[(u+v)\smile u']^\mathfrak{c}$ is onto to
$H^2(\Gamma,\mathfrak{c})$.
 This yields
$Q(u+v+w)=0$, with $Q$ being a submersion at $u+v+w$.


Therefore $Q^{-1}(0)$ is smooth of dimension $-\chi(\Gamma)\mathrm{dim}(G)$ near $u+v+w$. Since $v$, and consequently $w$, could be chosen arbitrarily small, this proves that $\bar{\mathcal{Q}}_{\phi}$ has a dense set of smooth points.
\end{pf}

\subsection{Proof of Theorem \ref{flexssimple}, flexibility part}
\label{proofflexsimpleflex}

Now we prove that \emph{Let $G$ be a semisimple real algebraic group. Let $\Gamma$ be the fundamental group of a closed surface of genus $\geq 2\mathrm{dim}(G)^2$. Let $\phi:\Gamma\to G$ be a homomorphism with reductive Zariski closure $H$. Assume that $\mathfrak{c}$, the center of the centralizer of $\phi(\Gamma)$, is balanced with respect to $\phi$. Then $\phi$ is flexible.
}

\begin{pf}
Proposition \ref{convex}
provides a dense set of smooth points for the central quadratic
model. Proposition \ref{extra} asserts that the full quadratic model has
the same property. Theorem \ref{dgms} transfers this property to a neighborhood of $\phi$ in $Hom(\Gamma,G)$. Thus $\phi$ is flexible.
\end{pf}

\section{Central extensions}
\label{ext}

In this section, we examine how smoothness lifts via central extensions. This is motivated by a technical difficulty in the proof of the remainder of Theorem \ref{flexssimple} (necessary condition for flexibility), see subsection \ref{proofflexsimplerig}. As a by product, this allows to generalize our results from semisimple to reductive ambient groups.

\subsection{Lifting homomorphisms}

\begin{lemma}
\label{centralextension}
Let $1\ra D\ra G\overset{\pi}{\ra} Q\ra 1$ be a central extension of Lie groups such that $D$ is discrete. Let $\Gamma$ be a finitely generated group. Let $\pi\circ:Hom(\Gamma,G)\ra Hom(\Gamma,Q)$, $\phi\mapsto\pi\circ\phi$. Then $\pi\circ$ is a local isomorphism of real analytic sets. In particular, it induces a local diffeomorphism $\mathcal{S}_{G}\to\mathcal{S}_{Q}$ (which needs not be onto).
\end{lemma}

\begin{pf}
Let us use the notation of subsection \ref{thevirtual}. Let $\langle S|N\rangle$ be a presentation of $\Gamma$. By noetherianity, near a point $\phi\in Hom(\Gamma,G)$ (resp. $\pi\circ\phi\in Hom(\Gamma,Q)$), a finite set $N'\subset N$ of relations suffice to define the real analytic subset $Hom(\Gamma,G)$ (resp. $Hom(\Gamma,Q)$). There is a commutative diagram
\begin{eqnarray*}
G^{S} &\stackrel{F_G}{\longrightarrow}& G^{N'}\\
\pi^{S}\downarrow&&\pi^{N}\downarrow\\
Q^{S} &\stackrel{F_Q}{\longrightarrow}& Q^{N'}
\end{eqnarray*}
in which the vertical arrows are local analytic isomorphisms. This shows that $\pi\circ:Hom(\Gamma,G)\ra Hom(\Gamma,Q)$ is a local analytic isomorphism. In particular, it maps smooth points to smooth points and conversely, and is a smooth diffeomorphism on the subset of smooth points.
\end{pf}

\begin{co}
\label{reductive}
Let $G$ be a connected, reductive real algebraic group with radical $R$ and $\pi:G\to G/R$. Let $\Gamma$ be a closed surface group of negative Euler characteristic. Let $\pi\circ:Hom(\Gamma,G)\ra Hom(\Gamma,G/R)$, $\phi\mapsto\pi\circ\phi$. Then $\pi\circ$ is an open map, its image is a union of connected components of $Hom(\Gamma,G/R)$. Furthermore, smooth homomorphisms $\Gamma\to G$ are exactly lifts of smooth homomorphisms $\Gamma\to G/R$, i.e. $\mathcal{S}_{G}=(\pi\circ)^{-1}(\mathcal{S}_{G/R})$.
\end{co}

\begin{pf}
Since $G$ is connected, reductive algebraic, $\mathfrak{r}$, the Lie algebra of $R$, is the center of $\mathfrak{g}$. Let $\mathfrak{g}=\mathfrak{r}\oplus\mathfrak{s}$ be a Levi decomposition of $\mathfrak{g}$, with $\mathfrak{s}$ semisimple. Then $[\mathfrak{g},\mathfrak{g}]=[\mathfrak{s},\mathfrak{s}]=\mathfrak{s}$, thus the Lie algebra $\mathfrak{r}\cap[\mathfrak{g},\mathfrak{g}]$ of $D=R\cap [G,G]$ is trivial, and $D$ is discrete.

Let $R'=R/D$ and $G'=G/D$. Then $R'\cap[G',G']=\{1\}$. Therefore the homomorphism $[G',G']\to G'/R'=G/R$ is injective. Its image has the same Lie algebra as $G/R$, it is an open subgroup of $G/R$, thus it is equal to $G/R$, since $G/R$ is connected. In other words, $[G',G']$ is isomorphic to $G/R$.

Since $R'$ is central in $G'$, the map $R'\times[G',G']\to G'$, $(r,h)\mapsto rh$, is a homomorphism, it is injective as $R'\cap[G',G']=\{1\}$, it is surjective as $[G',G']$ maps onto $G'/R'$, thus
\begin{eqnarray*}
G'\simeq R'\times[G',G'] \simeq R'\times G/R.
\end{eqnarray*}

The map $\pi:G\to G/R$ factors as $\pi''\circ\pi'$ where $\pi':G\to G'$ and $\pi'':G'\to G/R$. Since $D$ is discrete and $\pi'$ is a central extension, Lemma \ref{centralextension} applies, $\pi'\circ$ is locally homeomorphic, it maps smooth points of $Hom(\Gamma,G)$ to smooth points of $Hom(\Gamma,G')$ (and vice-versa) and is a smooth local diffeomorphism there.

Since $G'\simeq R'\times G/R$, $Hom(\Gamma,G')=Hom(\Gamma,R')\times Hom(\Gamma,G/R)$ and $\pi''\circ$ equals the projection onto the second factor. Since $Hom(\Gamma,R')$ is smooth, $\pi''\circ$ maps smooth points to smooth points and vice-versa. We conclude that $\mathcal{S}_{G}=(\pi\circ)^{-1}(\mathcal{S}_{G/R})$.
\end{pf}

\subsection{Virtual dimension for reductive groups}

Here, as announced in subsection \ref{virtual}, we justify the Definition \ref{defvirtual} of virtual dimensions for character varieties of reductive algebraic groups.

\begin{Prop}
\label{centralizercenter}
Let $G$ be a connected reductive real algebraic group with radical $R$. Let $\Gamma$ be a closed surface group. Let $\phi:\Gamma\ra G$ be a homomorphism such that the identity component of the centralizer of $\phi(\Gamma)$ is equal to $R$. Then $\phi$ belongs to the set $\mathcal{S}$ of smooth homomorphisms (see Definitions \ref{defrigid} and \ref{defvirtual}). In particular, Zariski dense homomorphisms belong to $\mathcal{S}$.
\end{Prop}

\begin{pf}
Let $\pi:G\ra G/R$. Under the assumption, the centralizer of $\pi\circ\phi(\Gamma)$ in $G/R$ is discrete, so $\pi\circ\phi$ is a smooth point of $Hom(\Gamma,G/R)$. According to Corollary \ref{reductive}, $\phi$ is a smooth point of $Hom(\Gamma,G)$ as well.
\end{pf}

\subsection{Proof of Theorem \ref{thm0}}
\label{proofthm1}

Here we complete the proof that {\em if $G$ is a connected reductive real algebraic group and  $\Gamma$ the fundamental group of a closed surface of genus $>1$, then the trivial homomorphism $\Gamma\to G$ can be deformed into flexible homomorphisms.}

\begin{pf}
Let $R$ denote the radical of $G$. According to Corollary \ref{cortrivial}, the trivial homomorphism belongs to the closure of $\mathcal{S}_{G/R}$ in $Hom(\Gamma,G/R)$. Corollary \ref{reductive} shows that its trivial lift to $G$ belongs to the closure of $\mathcal{S}_{G}$ in $Hom(\Gamma,G)$.
\end{pf}

\subsection{Flexibility versus rigidity for reductive representations}

\begin{Prop}
\label{clopen}
Let $G$ be a connected, reductive, real algebraic group. Let $\Gamma$ be
the fundamental group of a closed surface of genus $\geq
2\mathrm{dim}(G)^{2}$. Reductive homomorphisms $\Gamma\to G$ are either rigid or flexible.
\end{Prop}

\begin{pf}
Let $\mathcal{F}$ denote the set of flexible homomorphisms
$\Gamma\to G$. Recall that $\mathcal S$ denote the set of
homomorphisms which
are smooth points of $Hom(\Gamma,G)$. By definition,
$\mathcal{F}$ is open in $Hom(\Gamma,G)$ and contains $\mathcal S$
as a dense subset. So $\mathcal{F}$ is contained in the closure of
$\mathcal S$. Conversely, let $\phi$ be a reductive homomorphism
which belongs to the closure of $\mathcal S$. Let $R$ denote the
radical of $G$ and $\pi:G\ra G/R$. Then $\pi\circ\phi$ is reductive
again. Corollary \ref{reductive} tells that $Hom(\Gamma,G/R)$
admits smooth points in every neighborhood of $\pi\circ\phi$. Since
$G/R$ is semisimple, Theorem \ref{dgms} applies. The
Goldman-Millson-Simpson quadratic model at $\pi\circ\phi$ has at
least one smooth point. This is also valid for the central quadratic
model $\mathcal{Q}_{\pi\circ\phi}^{\mathfrak{c}}$. Proposition
\ref{extra} implies that the full quadratic model
$\mathcal{Q}_{\pi\circ\phi}$ has a dense set of smooth points. And
by Theorem \ref{dgms} again, smooth points are dense in a
neighborhood of $\pi\circ\phi$ in $Hom(\Gamma,G/R)$. By Corollary
\ref{reductive}, the same holds in a neighborhood of $\phi$
in $Hom(\Gamma,G)$, i.e. $\phi$ is flexible. This shows that
$RedHom(\Gamma,G)\cap\overline{\mathcal S}\subset \mathcal{F}$ where
$RedHom(\Gamma,G)$ is the set of reductive homomorphisms, which
implies that
$RedHom(\Gamma,G)\cap\mathcal{F}=RedHom(\Gamma,G)\cap\overline{\mathcal
S}$ is closed in $RedHom(\Gamma,G)$. Hence
$RedHom(\Gamma,G)\cap\mathcal{F}$ is both open and closed in
$RedHom(\Gamma,G)$, its complement consists of rigid homomorphisms only.
\end{pf}

\section{Amenable representations}
\label{amen}

Again, this section is motivated by a technical point in the proof of the remainder of Theorem \ref{flexssimple}. It will be proved that homomorphisms with reductive amenable Zariski closure are flexible. The generalization to nonreductive representations is treated in the next section.

\subsection{Bounded homomorphisms}

\begin{Prop}
\label{excompact}
Let $G$ be a compact semisimple real algebraic group. Let $\Gamma$ be a surface group. Let $\phi:\Gamma\to G$ be a homomorphism. Let $\mathfrak{c}$ be the center of the centralizer of $\phi$. Then $\mathfrak{c}$ is balanced with respect to $\phi$.
\end{Prop}

\begin{pf}
Every homomorphism to a compact group is reductive. Since $G$ is compact, the Killing form is definite, its restriction to all $(\mathfrak{g}_{\lambda,\br},J_{\lambda})$ is definite, so all Toledo invariants $T_{\lambda}$ vanish (Proposition \ref{abel}), the balance condition is automatically satisfied.
\end{pf}

\begin{co}
\label{compact} Let $G$ be a compact connected real algebraic group.
Let $\Gamma$ be the fundamental group of a closed surface of genus
$\geq 2\mathrm{dim}(G)^2$. Then every homomorphism $\Gamma\ra G$ is
flexible.
\end{co}

\begin{pf}
Let $R$ denote the center of $G$. View $G$ as a central extension of the compact semisimple real algebraic group $G/R$. Proposition \ref{excompact} allows to apply the flexibility part of Theorem \ref{flexssimple} (proved in \ref{proofflexsimpleflex}). This yields flexibility in $G/R$ and Corollary \ref{reductive} lifts it to $G$.
\end{pf}

\begin{co}
\label{bounded}
Let $G$ be a connected semisimple real algebraic group. Let $\Gamma$ be the fundamental group of a closed surface of genus $\geq 2\mathrm{dim}(G)^{2}$. Then every bounded homomorphism $\Gamma\ra G$ is flexible.
\end{co}

\begin{pf}
One can assume that $G$ is almost simple. Pick a maximal compact subgroup $K$ of $G$ containing $\phi(\Gamma)$. Since $G/K$ is contractible, $K$ is connected. According to Corollary \ref{compact}, $\phi$ can be deformed in $K$ to become Zariski dense in $K$. When the centralizer of $K$ in $\mathfrak{g}$ is trivial, $\phi$ belongs to $\mathcal{S}$. Otherwise $G/K$ is Hermitian-symmetric, the centralizer of $K$ in $\mathfrak{k}$ is the center of $\mathfrak{k}$, its image in the isotropy representation is generated by the complex structure $J$ on the tangent space at the point $[K]$. Its adjoint action on $\mathfrak{g}$ has only one non zero pair of roots, $\pm i$, one real root space $\mathfrak{p}=\mathfrak{k}^{\bot}$. Since the Killing form is definite on $\mathfrak{p}$, $J$ is tamed by $\Omega_{\lambda=i}$. The flat complex vectorbundle associated to the unitary representation $(\mathfrak{p},J)$ has vanishing first Chern class. Therefore the Toledo invariant $T_{i}$ vanishes, the balance condition is satisfied. Corollary \ref{abelian} implies flexibility in $G$.
\end{pf}

\subsection{Reductive amenable homomorphisms}

\begin{lemma}
\label{amenable}
Let $G$ be a real algebraic group. Let $\phi:\Gamma\to G$ be a homomorphism whose Zariski closure is reductive and amenable. Then the connected component of $\phi$ in the space $RedHom(\Gamma,G)$ of reductive homomorphisms $\Gamma\ra G$ contains a bounded homomorphism.
\end{lemma}

\begin{pf}
Let $H$ denote the Zariski closure of $\phi(\Gamma)$. Since $H$ is reductive, $H=SR$ where $R$, the radical, is connected and central, and $S$ is semisimple. The connected abelian algebraic group $R$ is isomorphic to a direct product $T\times V$ where $V$ is a vectorspace and $T$ is compact (\cite{Borel_LAG}, Proposition 8.15). $R\cap S$ is the center of $S$, thus is finite and contained in $T$. Therefore $TS\cap V=\{1\}$, and $H$ is isomorphic to $TS\times V$. Since $H$ is amenable, $S$ is compact, so is $TS$. For $(u,v)\in TS\times V$, let $\delta_{t}(u,v)=(u,tv)$. If $t\not=0$, this is an algebraic automorphism of $H$, thus the Zariski closure of $\phi_{t}=\delta_{t}\circ\phi$ is $\delta_{t}(H)=H$. The Zariski closure of $\phi_{0}$ is $TS$. In all cases, $\phi_{t}$ is reductive, so $t\mapsto\phi_{t}$ connects $\phi$ to the bounded homomorphism $\phi_{0}$ in $RedHom(TS\times V)\subset RedHom(\Gamma,G)$.
\end{pf}

\begin{co}
\label{coramenable}
Let $G$ be a connected reductive real algebraic group.
Let $\Gamma$ be the fundamental group of a closed surface of genus
$\geq 2\mathrm{dim}(G)^{2}$. Homomorphisms whose Zariski closure is
reductive and amenable are flexible.
\end{co}

\begin{pf}
Let $\phi:\Gamma\to G$ be a homomorphism whose Zariski closure is reductive and amenable.
 Lemma \ref{amenable} provides us with a bounded homomorphism $\phi_0$ in the connected component of $\phi$ in $RedHom(\Gamma,G)$.
 According to Corollary \ref{bounded}, $\phi_0$ is flexible. Proposition \ref{clopen} asserts that the
 set $\mathcal{F}\cap RedHom(\Gamma,G)$ of flexible reductive homomorphisms is open and closed, thus $\phi$ is flexible.
\end{pf}

This partial result will be generalized to arbitrary amenable homomorphisms in Corollary \ref{coramenablegeneral}.

\subsection{Proof of Theorem \ref{flexssimple}, rigidity part}
\label{proofflexsimplerig}

Here we complete the proof of Theorem \ref{flexssimple}, i.e. we prove that \emph{Let $G$ be a semisimple real algebraic group. Let $\Gamma$ be the fundamental group of a closed surface $\Sigma$ of genus $\geq 2\mathrm{dim}(G)^{2}$. Let $\phi:\Gamma\to G$ be a homomorphism whose Zariski closure $H$ is reductive. Assume that $\phi$ is flexible. Then the center of the centralizer of $\phi(\Gamma)$ is balanced with respect to $\phi$.} (The converse has been proven in subsection \ref{proofflexsimpleflex}).

\begin{pf}
Let $Z$ denote the centralizer of $H$, $Z^0$ the identity component of $Z$ and $\mathfrak{c}$ the center of $\mathfrak{z}$. We construct a small deformation $\phi'$ of $\phi$ whose centralizer in $\mathfrak{g}$ is exactly $\mathfrak{c}$. It suffices that $\phi'$ be Zariski dense in
$HZ^{0}$. Indeed, since $Z$ is reductive (Corollary \ref{centralizerisreductive}), $\mathfrak{c}$ is the center of $\mathfrak{z}$. Thus
\begin{eqnarray*}
\mathfrak{c}=Z_{\mathfrak{z}}(\mathfrak{z})=\mathfrak{z}\cap Z_{\mathfrak{g}}(\mathfrak{z})=Z_{\mathfrak{g}}(H)\cap Z_{\mathfrak{g}}(Z^0)=Z_{\mathfrak{g}}(HZ^0).
\end{eqnarray*}

Pick a Levi factor $L$ in $H$. Since the radical $R$ of $H$ is central and connected (\cite{Borel_LAG}, 11.21), it is contained in $Z^{0}$, so that $HZ^{0}=LZ^{0}$.

For clarity, let us treat first a simple case. Assume that $L\cap Z=\{1\}$ and $L$ and $Z$ are connected. Then $HZ$ is isomorphic to $L\times Z$. Let $pr_{Z}:LZ\ra Z$ and $pr_{L}:LZ\ra L$ denote the projections onto the factors. Since $\phi(\Gamma)$ is Zariski dense in $H$, $pr_{L}\circ\phi(\Gamma)$ is Zariski dense in $L$, and thus smooth as a homomorphism $\Gamma\to L$. The homomorphism $pr_{Z}\circ \phi$ has an abelian image. Corollary \ref{coramenable} implies that $pr_{Z}\circ \phi$ is flexible in $Z$, so it can be slightly deformed into a smooth homomorphism $\psi:\Gamma\ra Z$. Since $Hom(\Gamma,LZ)=Hom(\Gamma,L)\times Hom(\Gamma,Z)$, it follows that $(pr_{L}\circ\phi,\psi)$ is a smooth homomorphism $\Gamma\ra LZ$, which is close to $(pr_{L}\circ\phi,pr_{Z}\circ\phi)=\phi$. According to Corollary \ref{dense}, there is a Zariski dense $\phi':\Gamma\ra LZ$ close to $(pr_{L}\circ\phi,\psi)$ and thus to $\phi$.

In general, $F=L\cap Z^{0}$ is a finite group contained in the center of both $L$ and $Z^{0}$. The quotient $LZ^{0}/F$ is isomorphic to the direct product $(L/F)\times (Z^{0}/F)$. Let $pr_{L}:LZ^{0}\ra LZ^{0}/F\ra L/F$ and $pr_{Z}:LZ^{0}\ra LZ^{0}/F\ra Z^{0}/F$ denote composed projections. The image of the homomorphism $pr_{Z}\circ \phi : \Gamma\ra Z^{0}/F$ is Zariski dense in $pr_{Z}(R)$, thus is reductive and abelian. According to Corollary \ref{coramenable}, there exists a smooth homomorphism $\psi:\Gamma\ra Z^{0}/F$ nearby $pr_{Z}\circ \phi$. The image of $pr_{L}\circ \phi$ is Zariski dense in $L/F=H/R$. Then $pr_{L}\circ \phi$ and $\psi$ are smooth points of $Hom(\Gamma,L/F)$ and  $Hom(\Gamma,Z^{0}/F)$ respectively, thus $\bar{\phi}=(pr_{L}\circ \phi,\psi)$ is a smooth point of the product $Hom(\Gamma,L/F)\times Hom(\Gamma,Z^{0}/F)=Hom(\Gamma,L/F\times Z^{0}/F)$. Now $LZ^{0}$ is a finite central extension of $LZ^{0}/F$, thus, thanks to Lemma \ref{centralextension}, $\bar{\phi}$ lifts to a smooth homomorphism $\tilde{\phi}:\Gamma\to LZ^{0}$ close to $\phi$. The smallest open subgroup of $LZ^{0}$ containing $\phi(\Gamma)$ contains $Z^{0}$ (which is connected) and $H$, so this is $HZ^0=LZ^0$. Since $LZ^0$ has only finitely many open subgroups, $LZ^0$ is still the smallest open subgroup of $LZ^{0}$ containing $\tilde{\phi}(\Gamma)$. According to Proposition \ref{dense'}, there is a Zariski dense $\phi':\Gamma\ra LZ^{0}$ nearby $\tilde{\phi}$ and so nearby $\phi$. As we saw earlier, the centralizer of $\phi'$ equals $\mathfrak{c}$.

Since $\phi$ is flexible, so is $\phi'$. Corollary \ref{abelian} implies that $\mathfrak{c}$ is balanced with respect to $\phi'$. Since $\phi'$ is close to $\phi$, both homomorphisms give rise to the same Toledo invariants, $\mathfrak{c}$ is balanced with respect to $\phi$ as well.
\end{pf}

\section{Non reductive representations}
\label{nonreductive}

The fact that $\chigg=Hom(\Gamma,G)/G$ may be non Hausdorff greatly helps. It will allow us to continuously deform by conjugation non reductive representations to reductive ones.

\subsection{Closures of \texorpdfstring{$G$}{}-orbits in \texorpdfstring{$Hom(\Gamma,G)$}{}}

\begin{lemma}
\label{fixed}
Let $G$ be a semisimple real algebraic group, $H\subset G$ a non reductive
subgroup of $G$. There exists a one parameter subgroup $t\mapsto g_t$ of
$G$ such that the restriction of $Ad_{g_t}$ to $H$ converges as $t$ tends
to $+\infty$ to a homomorphism $T:H\to I$, where $I=Z_{G}(\{g_t\})$ is a product $I=G'K'A'$, where $G'$ is semisimple of non compact
type, $K'$ is compact, $A'$ is a real split torus, $G'$ commutes with $K'A'$
and $K'$ normalizes $A'$.
\end{lemma}

\begin{pf}
We rely on section 4.4 of \cite{Eberlein}. Let $X$ denote the
symmetric space of $G$, and let $X(\infty)$ denotes its visual
boundary, equipped with Tits' metric. By definition, the unipotent
radical $U$ of $H$ is non trivial. According to Proposition 4.4.3,
page 312 of \cite{Eberlein}, $U$ has fixed points on $X(\infty)$,
all contained in an open Tits ball of radius $\pi$. So does its normalizer
(Proposition 4.4.4 of \cite{Eberlein}), so $H$ has a fixed point $x$
in $X(\infty)$. Next we use section 2.17 of \cite{Eberlein}. Pick an
origin $o$ in $X$. The unit vector $u$ at $o$ pointing towards $x$
corresponds to an element, still denoted by $u$, of the
$\mathfrak{p}$ part of the Cartan decomposition of $\mathfrak{g}$ at
$o$. According to Proposition 2.17.5, page 104 of \cite{Eberlein},
for every $g\in H$, the limit
\begin{eqnarray*}
T_{x}(g)=\lim_{t\ra +\infty} e^{-tu}ge^{tu}
\end{eqnarray*}
exists, and defines a group homomorphism $T_{x}:H\to I =Z_{G}(u)$.
Furthermore, $I=K_{x}A_{x}$ where $K_{x}=I\cap K$ is compact
and $A_{x}=exp(Z_{\mathfrak{g}}(u)\cap\mathfrak{p})$, so $I$ is
stable under the Cartan involution corresponding to $o$. It follows
that the orbit $Y$ of $o$ under $I$ is a totally geodesic subspace
of $X$. Let $F$ denote the Euclidean factor of the de Rham
decomposition $Y=Y'\times F$ of $Y$. The isometry group $G'$ of $Y'$
being generated by geodesic symmetries, embeds in $G$. For the same
reason, the group $A'$ of translations of $F$ embeds in $G$. $G'$
and $A'$ are contained in $I$. Since $G'\times A'$ acts
transitively on $Y=Y'\times F$, it is cocompact in $I$. The
kernel of the isometric action of $I$ on $Y'$ is a subgroup of
$Isom(F)$ which contains $A'$, i.e. a semi-direct product $K'\ltimes
A'$ with $K'$ compact. Therefore $I$ is isomorphic to $G'\times
(K'\ltimes A')$.
\end{pf}

The next lemma is a variant of Lemma \ref{fixed}, adjusted for an induction procedure.

\begin{lemma}
\label{cascade} Let $G$ be a semisimple real algebraic group, let
$G_1$ and $G_2$ be commuting  subgroups of $G$, such that $G_1$ is
semisimple of non compact type, and $G_2$ contains no unipotent
elements. Let $H\subset G_1 G_2$ be a non reductive subgroup of $G_1
G_2$. There exists a one parameter subgroup $t\mapsto g_t$ of $G_1$
such that the restriction of $Ad_{g_t}$ to $H$ converges as $t$
tends to $+\infty$ to a homomorphism $T:H\to G'_1 G'_2$, where
$G'_1$ is a proper subgroup of $G_1$ and is semisimple of non
compact type, $G'_2$ commutes with $G'_1$ and contains
no unipotent elements.
\end{lemma}

\begin{pf}
Let $U$ be the unipotent radical of $H$. Its projection to $G_2$ is
unipotent, thus trivial,  so $U\subset H\cap G_1$. In particular,
the projection $H_1$ of $H$ to $G_1$ is not reductive. Lemma
\ref{fixed} provides the needed one-parameter group $t\mapsto g_t$.
$Ad_{g_t}$ being trivial on $G_2$, its restriction to $H_1 G_2
\supset H$ converges to a homomorphism $T:H_1 G_2 \to
Z_{G_1}(\{g_t\})G_2$. Furthermore, $Z_{G_1}(\{g_t\})=G'_1 K'A'$
where $G'_1$ is a proper subgroup of $G_1$ and is semisimple of non
compact type, and $K'A'$ commutes with $G'_1$ and contains no
unipotent elements. Since $K'A'$ commutes with $G_2$, $G'_2
=K'A'G_2$ does not contain any unipotent elements either.
\end{pf}

\begin{Prop}
\label{squeeze}
Let $G$ be a semisimple real algebraic group. Let $\Gamma$ be a finitely generated group. Let $G$ act on $Hom(\Gamma,G)$ by conjugation. Then every orbit closure contains at least one reductive homomorphism. More precisely, given a homomorphism $\phi:\Gamma\to G$ with Zariski closure $H$, there exists a homomorphism $T:H\to G$ such that $T\circ\phi$ is reductive and belongs to the closure of the $G$-orbit of $\phi$.
\end{Prop}

\begin{pf}
Let $H$ be the Zariski closure of $\phi(\Gamma)$. If $H$ is not reductive, apply Lemma \ref{cascade} to $G_1 =G$ and trivial $G_2$, and get $T_{1}:H\to G'_1 G'_2$. If $T_{1}(H)$ is not reductive, replace $\phi=\phi_0$ with $\phi_1 =T_{1}\circ\phi$ and apply Lemma \ref{cascade} again, and so on. When the procedure stops, we get a reductive homomorphism $\psi=\phi_k =T_{k}\circ\cdots\circ T_{1}\circ\phi$. By construction, each $\phi_j$ belongs to the closure of the $G$-orbit of $\phi_{j-1}$. Thus $\psi$ belongs to the closure of the $G$-orbit of $\phi$.
\end{pf}

\begin{co}
\label{coramenablegeneral} Let $G$ be a connected reductive real algebraic group. Let $\Gamma$ be the fundamental  group of a closed surface of genus
$\geq 2\mathrm{dim}(G)^{2}$. Homomorphisms whose Zariski closure is
amenable are flexible.
\end{co}

\begin{pf}
This is a combination of Proposition \ref{squeeze}, openness of flexibility and Corollary \ref{coramenable}. Indeed, images of amenable groups are amenable.
\end{pf}

\subsection{Proof of Theorem \ref{thm1}}

{\em Let $G$ be a connected reductive real algebraic group. Let $\Gamma$ be the fundamental group of a closed surface $\Sigma$ of genus $\geq 2\mathrm{dim}(G)^2$. Homomorphisms $\Gamma\to G$ are either flexible or rigid.}

\begin{pf}
Assume first that $G$ is semisimple. The sets $\mathcal{F}$ of flexible homomorphisms and $\mathcal{R}$ of rigid homomorphisms are both open and $G$-invariant. The complement of $\mathcal{F}\cup\mathcal{R}$ in $Hom(\Gamma,G)$ is closed and $G$-invariant. Assume that it is non empty. According to Proposition \ref{squeeze}, the closure of the orbit of any of its elements contains a reductive homomorphism, which, according to Proposition \ref{clopen}, belongs to $\mathcal{F}\cup\mathcal{R}$, contradiction. We conclude that $Hom(\Gamma,G)=\mathcal{F}\cup\mathcal{R}$.

If $G$ is merely reductive, let $R$ denote its radical, $\pi:G\ra
G/R$ and denote by $\pi\circ : Hom(\Gamma,G)\ra Hom(\Gamma,G/R)$ the
induced map. From Corollary \ref{reductive}, we know that smooth
points correspond under $\pi\circ$, i.e.
$\mathcal{S}_{G}=(\pi\circ)^{-1}\mathcal{S}_{G/R}$. Since $\pi\circ$
is an open map, it commutes with closures and interiors, so
$\mathcal{F}_{G}=(\pi\circ)^{-1}\mathcal{F}_{G/R}$ and
$\mathcal{R}_{G}=(\pi\circ)^{-1}\mathcal{R}_{G/R}$.  Since $G/R$ is
semisimple, we know that
$Hom(\Gamma,G/R)=\mathcal{F}_{G/R}\cup\mathcal{R}_{G/R}$. It follows
that $Hom(\Gamma,G)=\mathcal{F}_{G}\cup\mathcal{R}_{G}$.
\end{pf}

\subsection{Proof of Corollary \ref{rigidtt}}
\label{proofrigidtt}

The following result is a special case of Theorem 5 of \cite{BIW}.

\begin{thm}
\label{cortight}
Let $\Gamma$ be surface group, $\phi:\Gamma\to G$ a homomorphism. Let $G$ act symplectically on $(V,\Omega)$. Let $\rho:\Gamma\to Sp(V,\Omega)$ be a homomorphism. Assume equality holds in the Milnor-Wood inequality, i.e.
\begin{eqnarray*}
|T_{\rho}|=-\chi(\Gamma)\mathrm{dim}(V).
\end{eqnarray*}
Then the Zariski closure of $\rho(\Gamma)$ in $Sp(V,\Omega)$ is reductive, the symmetric space associated to its semisimple Levi factor is Hermitian of tube type.
\end{thm}

Using this theorem, we can prove Corollary \ref{rigidtt}, i.e. \emph{Let $G$ be a connected reductive real algebraic group. Let $\Gamma$ be the fundamental group of a closed surface of genus $\geq 2\mathrm{dim}(G)^2$. Let $\phi:\Gamma\to G$ be a homomorphism with Zariski closure $H$. If $\phi$ is not flexible, then $H$ acts transitively on a tube type Hermitian symmetric space, producing a maximal representation of $\Gamma$.}

\begin{pf}
Assume first that $G$ is semisimple and $H$ is reductive. The centralizer
$\mathfrak{z}$ of $H$ in $\mathfrak{g}$ is reductive (Corollary \ref{centralizerisreductive}). If $\phi$ is not
flexible, then $\mathfrak{z}$ cannot be semisimple (Proposition
\ref{surface1}) so it has a nontrivial center $\mathfrak{c}$. Furthermore
(sufficient condition for flexibility in Theorem \ref{flexssimple}, proved
in subsection \ref{proofflexsimpleflex}) there is at least one nonzero root $\lambda$ of $\mathfrak{c}$ on $\mathfrak{g}$ for which
\begin{eqnarray}
4|c_1 (E_{\lambda})|= -\chi(\Gamma)\mathrm{rank}(E_{\lambda}).\label{MW2}
\end{eqnarray}
Here $\mathfrak{g}_{\lambda,\br}$ is the real root space corresponding to
$\lambda$, $\Omega_{\lambda}$ the symplectic form on it, $E_{\lambda}$ the
symplectic bundle on $\Sigma$ induced by an equivariant map
$\tilde{\Sigma}\to X(\mathfrak{g}_{\lambda,\br},\Omega_{\lambda})$. Since
$H$ centralizes $\mathfrak{c}$, it preserves the splitting of
$\mathfrak{g}$ and acts symplectically on $\mathfrak{g}_{\lambda,\br}$. Its
image $I$ in $Sp(\mathfrak{g}_{\lambda,\br},\Omega_{\lambda})$ coincides
with the Zariski closure of the image of $\Gamma$ in
$Sp(\mathfrak{g}_{\lambda,\br},\Omega_{\lambda})$, which, according to
Theorem \ref{cortight}, is reductive with a Levi factor which is Hermitian
of tube type. $I$ is a quotient of $H$, thus its Levi factor is a product
of simple factors of a semisimple Levi factor of $H$. Thus at least one of them gives rise to a Hermitian symmetric space of tube type, with a maximal representation of $\Gamma$, thanks to equation (\ref{MW2}).

If $H$ is not reductive, Proposition \ref{squeeze} yields a homomorphism $T:H\to G$ such that $\phi'=T\circ\phi$ is a limit of conjugates of $\phi$ and has a reductive Zariski closure $H'$. If $\phi$ is not flexible, neither is $\phi'$, so one of the simple factors of the semisimple Levi factor of $H'$ has a tube type Hermitian symmetric space. $H$, which surjects onto $H'$, acts transitively on this symmetric space.

In general, let $G$ be connected reductive with radical $R$, and $\pi:G\to G/R$. Then $G/R$ is semisimple. According to Corollary \ref{reductive}, the projected homomorphism $\pi\circ\phi$ is not flexible. Thus $\pi(H)$ acts transitively on a tube type Hermitian symmetric space, so does $H$.
\end{pf}

\section{Surface groups in real rank 1 simple Lie groups}
\label{one}

Corollary \ref{rankone} is proven here.

\subsection{Consequences of general results}

Let $G$ be a simple Lie group of real rank 1. Let $\Gamma$ be the fundamental group of a closed surface of large genus $\geq 2\mathrm{dim}(G)^2$. Let $\phi:\Gamma\to G$ be a homomorphism.

If the Zariski closure of $\phi(\Gamma)$ is not reductive, then it is contained in a proper parabolic subgroup. In rank 1, there is only one, the minimal parabolic subgroup $AN$, which is amenable. Thus Corollary \ref{coramenable} applies, and $\phi$ is flexible.

Assume that the Zariski closure $H$ of $\phi(\Gamma)$ is reductive nonamenable. According to Lemma \ref{amencentr}, the centralizer of a Levi factor of $H$ is compact, so is the centralizer $C$ of $H$.

Assume $\phi$ is not flexible in $G$. Let $S$ denote the product of noncompact factors of a Levi factor of $H$. The symmetric space $Y$ of $S$ is totally geodesic in $X$, it has rank $1$, so $S$ has only one factor. According to Corollary \ref{rigidtt}, $Y$ has to be Hermitian of tube type. It must be a complex hyperbolic space. Only the complex hyperbolic line is of tube type. It follows that $S$ is locally isomorphic to $SL(2,\br)$. Totally geodesic 2-planes in $X=H^{m}_{\bh}$ are of the following 2 types
\begin{enumerate}
  \item $\bc$-line, i.e. subsets of $\bF$-lines, $\bF=\bc$, $\bh$, $\bo$.
  \item $\br$-planes, spanned by orthonormal tangent vectors $a$, $b$ with $b\bot \bF a$.
\end{enumerate}

In each case, we shall compute the center of the centralizer. It turns out to always be $0$ or $1$-dimensional. A trivial center implies flexibility (Proposition \ref{surface1}). According to Theorem \ref{flexssimple}, when the center of the centralizer is $1$-dimensional, a necessary and sufficient condition for flexibility is that the direct sum of symplectic representations $\rho_\lambda$ be non maximal.

In the computation of the Toledo invariant of $\rho_\lambda$, only the $S$-component of $\phi$ plays a role. Indeed, since $Z_{G}(S)$ is compact, one can always choose an equivariant map $\tilde{\Sigma}\ra Sp(\mathfrak{g}_{\lambda,\br})$ which factorizes via the symmetric space of $S$.

\subsection{Notation}

For $\bF=\br$, $\bc$, $\bh$, $\bo$, let $O(m,p,\bF)$ denote
respectively $O(m,p)$, $U(m,p)$, $Sp(m,p)$,  $F_{4}^{-20}$ with a
corresponding notation for Lie algebras.

For instance, a nonamenable Lie group $H$ that stabilizes a real plane contains
$S=SO(2,1,\br)\subset SO(2,1,\bF)\subset SO(m,1,\bF)$. If $H$ stabilizes a
complex line, it contains $S=SU(1,1)=SO(1,1,\bc)\subset SO(1,1,\bF)\subset SO(m,1,\bF)$.

The following lemma will be useful.

\begin{lemma}
\label{splitso}
For $\bF=\br$, $\bc$, $\bh$, for $0<p<m$,
\begin{eqnarray*}
\mathfrak{o}(m,q,\bF)=\mathfrak{o}(m-p,\bF)\oplus \mathfrak{o}(p,q,\bF)\oplus Hom_{\bF}(\bF^{m-p} ,\bF^{p+q}).
\end{eqnarray*}

\end{lemma}

\begin{pf}
Let $e=\begin{pmatrix}
I_p & 0 \\
0 & -I_q
\end{pmatrix}$. Then a $(m-p)\oplus (p+q)$-block matrix $M=\begin{pmatrix}
A & B \\
C & D
\end{pmatrix}$ with entries in $\bF$ belongs to $\mathfrak{o}(m,q,\bF)$ iff $A\in \mathfrak{o}(m-p,\bF)$, $D\in \mathfrak{o}(p,q,\bF)$ and $B=-C^{*}e=-\bar{C}^{\top}e$. An element $M=\begin{pmatrix}
A & 0 \\
0 & D
\end{pmatrix}\in O(m-p,\bF)\times O(p,q,\bF)$ acts on $N=\begin{pmatrix}
0 & -C^{*}e \\
C & 0
\end{pmatrix}$ by
$$
Ad_{M}(N)=\begin{pmatrix}
0 & -A^{-1}C^{*}eD \\
D^{-1}CA & 0
\end{pmatrix},
$$
justifying the notation for the complement of $\mathfrak{o}(m-p,\bF)\oplus \mathfrak{o}(p,q,\bF)$ in $\mathfrak{o}(m,1,\bF)$. For future use, we record the following formula. The Lie bracket of two elements $\begin{pmatrix}
0 & -B^{*}e \\
B & 0
\end{pmatrix}$ and $\begin{pmatrix}
0 & -C^{*}e \\
C & 0
\end{pmatrix}$ of $Hom_{\bF}(\bF^{m-p} ,\bF^{p+q})$ equals
$
\begin{pmatrix}
C^{*}eB-B^* eC  &  0 \\
0   & CB^* e-BC^* e
\end{pmatrix}$.
\end{pf}

\subsection{Homomorphisms that stabilize an \texorpdfstring{$\br$}{}
-plane}

Let $\bF=\br$, $\bc$ or $\bh$. The centralizer of $SO(2,1,\br)$ in $O(2,1,\bF)$ is the group of multiples of identity by units of $\bF$, $O(1,\bF)$. The
centralizer of $SO(2,1,\br)$ in $O(m,1,\bF)$ stabilizes a $\bF$-plane,
thus it is equal to $O(1,\bF)\times O(m-2,\bF)$.

Case $\bF=\bo$. Recall that the maximal compact subgroup $K$ of
$G=F_{4}^{-20}$ is isomorphic to $Spin(9)$, and the action of $K$ on the
tangent space $T$ to $G/K$ is the 16-dimensional spin representation of
$Spin(9)$. Let $v$, $w$ be orthonormal vectors in $T$ with $w\in(\bo
v)^{\bot}$. The fixator of $v$ in this representation equals $M=Spin(7)$,
$T$ splits as $T=\br v \oplus \Im m(\bo)v \oplus (\bo v)^{\bot}$, on which
$M$ acts respectively by the trivial, standard $SO(7)$ and spin
representations. The fixator of $v$ and $w$ equals the fixator of $w$ in
the spin representation of $Spin(7)$, i.e. $G_2$.

To sum up, the centralizer $Z$ of $S=SO(2,1)$ equals
\begin{itemize}
  \item $O(m-2)$ if $\bF=\br$,
  \item $S(U(1)\times U(m-2))$ if $\bF=\bc$,
  \item $Sp(1)\times Sp(m-2)$ if $\bF=\bh$,
  \item $G_2$ if $\bF=\bo$.
\end{itemize}
Then, the center of $\mathfrak{z}$ is
\begin{itemize}
  \item $\mathfrak{so}(2)$ if $\bF=\br$ and $m=4$,
  \item trivial if $\bF=\br$ and $m\not=4$,
  \item trivial if $\bF=\bc$ and $m=2$,
  \item $\mathfrak{u}(1)$ (the center of $S(U(1)\times U(m-2))$), if $\bF=\bc$, $m\geq 3$,
  \item trivial if $\bF=\bh$,
  \item trivial if $\bF=\bo$.
\end{itemize}

Thus, if genus is large enough, $\phi$ is flexible, except possibly in $SO(4,1)$ or in $SU(m,1)$, $m\geq 3$. It has been known for a long time that flexibility holds in theses cases (these homomorphisms can be deformed by bending, see \cite{Th}, \cite{Apa}), but we include a proof along the lines of Theorem \ref{flexssimple}, for completeness' sake.

\subsubsection{4-dimensional real hyperbolic space}
\label{4real}
The splitting of $\mathfrak{so}(4,1)$ in root spaces according to the centralizer $\mathfrak{so}(2)$ is given by Lemma \ref{splitso}. There is only one space $\mathfrak{g}_{\lambda,\br}$, it is $Hom(\br^2 ,\br^3)$. Represent elements $C\in Hom(\br^2 ,\br^3)$ as pairs of vectors $(c_1 ,c_2)$ of $\br^3$. Then the $\mathfrak{so}(2)$-component of the Lie bracket $[B,C]$ is $((b_1 ,b_2),(c_1 ,c_2))\mapsto c_1 \cdot b_2 -c_2 \cdot b_1$ (here, $\cdot$ denotes the $SO(2,1)$-invariant symmetric bilinear form on $\br^3$). This is a symplectic structure. The subspaces $R=\{c_2 =0\}$ and $R^* =\{c_1 =0\}$ are Lagrangian and $SO(2,1)$-invariant. Lemma \ref{toledochern} implies that the Toledo invariant attached to this symplectic structure vanishes. Therefore flexibility holds.

\subsubsection{Real planes in complex hyperbolic space}

Again, Lemma \ref{splitso} provides us with the only real root space of the centralizer, $Hom_{\bc}(\bc^{m-2} ,\bc^3)$. Let $Z=\begin{pmatrix}
\frac{3i}{m+1}I_{m-2} & 0   \\
0    &  \frac{(2-m)i}{m+1}I_{3}
\end{pmatrix}$ be a basis of the centralizer. Then $ad_{Z}$ is nothing but the obvious complex structure on $Hom(\bc^{m-2} ,\bc^3)$. Then $(B,C)$ $\mapsto$ $\Omega_{Z}(ad_{Z}B,C)$ is the direct sum of $m-2$ copies of the $SU(2,1)$-invariant Hermitian form on $\bc^3$. Let $R=Hom_{\br}(\br^{m-2} ,\br^3)\otimes\bc$ $\subset$ $Hom_{\bc}(\bc^{m-2} ,\bc^3)$ denote the subspace of real matrices, and $R^* =iR$. These are two Lagrangian, $SO(2,1)$-invariant subspaces. Lemma \ref{toledochern} implies that the Toledo invariant attached to this symplectic structure vanishes.

\subsection{Homomorphisms that stabilize a \texorpdfstring{$\bc$}{}
-line}

Let $\bF=\bc$ or $\bh$. The centralizer of $SU(1,1)$ in $SO(1,1,\bh)=Sp(1,1)$ equals $U(1)$ (the centralizer of $SU(1,1)$ in $U(1,1)$). The centralizer of $SU(1,1)$ in $SO(m,1,\bF)$ fixes a $\bc$-line pointwise, thus it stabilizes the $\bF$-line that contains it. Thus it is equal to $U(m-1)$ (complex case), $U(1)\times Sp(m-1)$ (quaternionic case).

Case $\bF=\bo$. Let $T$ denote the tangent space to octonionic hyperbolic
plane. Let $v$, $w$ be orthonormal vectors in $T$ with $w\in\bo v$. The
fixator of $v$ in the isotropy representation equals $Spin(7)$, $T$ splits as $T=\br v \oplus \Im m(\bo)v \oplus (\bo v)^{\bot}$, on which $M$ acts respectively by the trivial, standard $SO(7)$ and spin representations. The fixator of $v$ and $w$ equals the fixator of $w$ in the standard representation of $SO(7)$, i.e. a double cover of $SO(6)$, $Spin(6)$. Therefore the centralizer of $SU(1,1)$ in $F_{4}^{-20}$ is isomorphic to $Spin(6)$.

To sum up, the centralizer $Z$ of $S=SU(1,1)$ is
\begin{enumerate}
  \item $S(U(1)\times U(m-1))$ if $\bF=\bc$,
  \item $U(1)\times Sp(m-1)$ if $\bF=\bh$,
  \item $Spin(6)$ if $\bF=\bo$.
\end{enumerate}
Then, the center of $\mathfrak{z}$ is
\begin{enumerate}
  \item $\mathfrak{u}(1)$ (the center of $S(U(1,1)\times U(m-1))$) if $\bF=\bc$,
  \item $\mathfrak{u}(1)$ (the center of $U(1,1)\subset Sp(1,1)\subset Sp(m,1)$) if $\bF=\bh$,
  \item trivial if $\bF=\bo$.
\end{enumerate}

\subsubsection{Octonionic case}

A trivial center implies flexibility.

\subsubsection{Quaternionic case}

Choose $Z=\begin{pmatrix}
0 & 0 \\
0 & iI_2
\end{pmatrix}$ as a basis of the centralizer. Once more, Lemma \ref{splitso} provides us with the unique real root space, $Hom(\bF^{m-1},\bF^{2})$, on which $ad_{Z}$ acts by multiplication by $i$, i.e. as the obvious complex structure.

We continue with $\bF=\bh$. The Killing form on $\mathfrak{sp}(m,1)$ is proportional to $M\mapsto \Re e(\mathrm{Trace}_{\bh}(M^{2}))$.
For $X=\begin{pmatrix}
0 & -B^{*}e \\
B & 0
\end{pmatrix}$ and $Y=\begin{pmatrix}
0 & -C^{*}e \\
C & 0
\end{pmatrix}$, up to a constant,
\begin{eqnarray*}
\Omega_Z (X,Y)&=&Z\cdot[X,Y]\\
&=&\Re e(\mathrm{Trace}_{\bh}(Z[X,Y]))\\
&=&\Re e(i\,\mathrm{Trace}_{\bh}(CB^* e-BC^* e))\\
&=&2\Re e(i\,\mathrm{Trace}_{\bh}(CeB^*)).
\end{eqnarray*}
This alternating form is, in an $Sp(1,1)$-invariant manner, the direct sum of $m-1$ copies of a symplectic structure $\Omega$ associated with a $\bc$-Hermitian form on $\bh^{2}$ (it turns out to be the $Sp(1,1)$-invariant $\bh$-Hermitian form). Therefore we need only compute the Toledo invariant in case $m=2$. In this case, $Sp(1,1)=SO(4,1)$, $SU(1,1)$ is a diagonal $SO(2,1)$. The calculation, done in paragraph \ref{4real}, shows that the Toledo invariant vanishes. Flexibility follows. This example, studied in \cite{KPu}, is the starting point of the present work.

\subsubsection{Complex case}

Now $\bF=\bc$. Since the Killing form on $\mathfrak{su}(m,1)$ is proportional to $M\mapsto \Re e(\mathrm{Trace}_{\bc}(M^{2}))$, for $X=\begin{pmatrix}
0 & -B^{*}e \\
B & 0
\end{pmatrix}$ and $Y=\begin{pmatrix}
0 & -C^{*}e \\
C & 0
\end{pmatrix}$,
\begin{eqnarray*}
\Omega_Z (X,Y)&=&Z\cdot[X,Y]=\Re e(\mathrm{Trace}_{\bc}(Z[X,Y]))\\
&=&\Re e(i\,\mathrm{Trace}_{\bc}(CB^* e-BC^* e))=2\Re e(i\,\mathrm{Trace}_{\bc}(B^* eC)).
\end{eqnarray*}
This alternating form is, in an $SU(1,1)$-invariant manner, the direct sum of $m-1$ copies of the symplectic structure $\Omega$ associated with the $SU(1,1)$-invariant Hermitian form on $\bc^{2}$. Therefore the corresponding Toledo invariant is equal to $m-1$ times the Toledo invariant computed in Example \ref{su(1,1)->su(2,1)}. Thus it is maximal if and only if the projection $\phi'$ of $\phi$ to the Levi factor $S=SU(1,1)$ of its Zariski closure is maximal. According to Goldman's theorem \ref{PU(1,1)}, this happens if and only if $\phi'$ is discrete and cocompact. We conclude that $\phi$ is rigid if and only if $\phi(\Gamma)$ is discrete and cocompact in $S(U(1,1)\times U(m-1))$.

\subsection{Proof of Corollary \ref{rankone}}

It follows from Theorem \ref{flexssimple} and the case by case study of the preceding subsections.

\begin{rem}
Real rank one allows many shortcuts in the proof of Theorem \ref{flexssimple}.
\end{rem}

\section{Appendix 1 : Computing the cup-product}
\label{cupproduct}

\subsection{W. Meyer's signature formula}

Given a linear form $\ell$ on the second cohomology, the cup-product composed with $\ell$ is a real valued quadratic form. Its signature will be given by the following lemma, which is a special case of W. Meyer's signature formula, \cite{Meyer}).

\begin{Def}
\label{deftame} The {\em first Chern class} of a symplectic vector bundle
$(E,\Omega)$ is defined as follows : pick a smooth complex structure $J$ on the fibres
which is tamed by $\Omega$, i.e. such that the quadratic form $e\mapsto
\Omega(e,Je)$ is positive definite. Up to homotopy, $J$ is unique.
Therefore the first Chern class of the complex vector bundle $(E,J)$
does not depend on the choice of $J$.
\end{Def}

\begin{lemma}
\label{index} {\em (\cite{Meyer})}.
Let $(E,\Omega)$ be a flat symplectic vector bundle over a closed
oriented $2$-manifold $\Sigma$. The quadratic form
$\kappa(a,b)=\int_{S}\Omega(a\smile b)$ on $H^{1}(\Sigma,E)$ is
nondegenerate of signature $4c_1 (E,\Omega)$.
\end{lemma}

\begin{pf}
Since \cite{Meyer} is not so easily available, we include a short
proof of the special case we need. See also \cite{Py}.

In general, if $E$ is a flat vector bundle over a closed
$n$-manifold $X$, and $\Omega$ is a nondegenerate parallel bilinear
form on $E$, Poincar\'e duality ensures that the pairing
$(a,b)\ra\int_{X}\Omega(a\smile b)$ is nondegenerate on
$H^{k}(X,E)\otimes H^{n-k}(X,E)$. This can be seen as follows. Pick
a smooth Euclidean structure $G$ on the fibres of $E$. View $G$ as a
section of $Hom(E,E^*)$ since an inner product on a vector space
determines an isomorphism to its dual. This defines an $L^2$ norm
$\int_{X}\langle G\alpha,\alpha\rangle \,vol_X$ on $E$-valued forms
on $X$. With respect to this Hilbert space structure, the adjoint of
the exterior differential is
\begin{eqnarray*}
d^{*}=\pm G^{-1}\star d\star G,
\end{eqnarray*}
where in this equation, $d$ denotes exterior differentiation of
$E$-valued forms, and $\star$ denotes the Hodge $\star$ operator,
trivially extended from $\Lambda^* T^*X$ to $\Lambda^* T^*X \otimes
E$. View $\Omega$ as a parallel section of $Hom(E,E^*)$. Then
$\Omega$ commutes with exterior differentiation. Therefore
\begin{eqnarray*}
d^{*}=\pm G^{-1}\Omega \star d\star \Omega^{-1} G=\pm
\star_{\Omega}d\star_{\Omega},
\end{eqnarray*}
where $\star_{\Omega}$ maps harmonic forms to harmonic forms, and
satisfies
\begin{eqnarray*}
\Omega(\alpha\wedge\star_{\Omega}\alpha)=|\alpha|^2 vol_X
\end{eqnarray*}
pointwise. The integrated version
\begin{eqnarray*}
\int_{X}\Omega(\alpha\wedge\star_{\Omega}\alpha)=||\alpha||_{L^2}^2
\end{eqnarray*}
implies that $\kappa$ is nondegenerate on harmonic forms, and thus
on cohomology.

Assume $X=\Sigma$ is $2$-dimensional and $\Omega$ is skew-symmetric.
Let $J$ be a smooth complex structure on the fibres of $E$ tamed by
$\Omega$. Choose $G(e)=\Omega(e,Je)$. Then
$$|\alpha(v)|^2_G vol(v,\star v)=\Omega(\alpha(v),
(\star_\Omega\alpha)(v))=\Omega(\alpha(v),J\alpha(v))vol(v,\star
v).$$ Hence $\star_{\Omega} =\star\otimes J$ and it is an involution
on $E$-valued $1$-forms. Let $\pi_+:\Lambda^1 T^*\Sigma \otimes E
\ra \Lambda^+$ denote projection onto the $+1$-eigenspace of
$\star_{\Omega}$, and $d^+ =\pi_+ \circ d$. If an $E$-valued 0-form
$e$ satisfies $d^- e=0$, then $de=d^+e+d^-e=d^+e=\star_\Omega
d^+e=\star_{\Omega}de$ is harmonic (for
$(\star_{\Omega}d\star_{\Omega}d+d\star_{\Omega}d\star_{\Omega})(de)=0$)
and exact, thus vanishes. In other words, $\mathrm{ker}d^-
=\mathrm{ker}d\simeq H^0 (\Sigma,E)$. Let $\alpha$ be an
antiselfdual $E$-valued 1-form, i.e., $\star_{\Omega}\alpha=-\alpha$
and $e$ an $E$-valued 0-form. Then $\langle \alpha,d^-
e\rangle_{L^2}=\langle \alpha,de\rangle_{L^2}$. For
$$\langle \alpha, de\rangle=\langle -\star_{\Omega}\alpha,
d^+e+d^-e\rangle=\langle -\alpha,
\star_{\Omega}d^+e+\star_{\Omega}d^-e\rangle=\langle
-\alpha,d^+e-d^-e\rangle,$$ which implies $\langle
\alpha,d^+e\rangle=0$.

 Thus if $\alpha$ is
orthogonal to the image of $d^-$, then
$\alpha=-\star_{\Omega}\alpha$ is coclosed (for $\langle
\alpha,de\rangle=\langle \delta \alpha,e\rangle=0$ for all $e$, and
hence $\delta\alpha=0$) and also closed (for
$d\alpha=-d\star_{\Omega}\alpha=-\star_{\Omega}\star_{\Omega}d\star_{\Omega}\alpha=-\star_{\Omega}\delta\alpha=0$),
thus harmonic. In other words, the cokernel of $d^-$ is the space
$\mathcal{H}^-$ of antiselfdual harmonic 1-forms. Therefore
$\mathrm{Index}(d^-)=\mathrm{dim}ker-\mathrm{dim}coker=\mathrm{dim}H^0
(\Sigma,E)-\mathrm{dim}\mathcal{H}^-$. Since $d^-$ is the
$\bar{\partial}$ operator of the complex vector bundle $(E,J)$,
\begin{eqnarray*}
\mathrm{Index}(d^-)=\chi(\Sigma)\mathrm{rank}_{\bc}(E,J)+2c_1 (E,J),
\end{eqnarray*}
(see for instance \cite{MDS}). Similarly,
\begin{eqnarray*}
\mathrm{Index}(d^+)&=&\chi(\Sigma)\mathrm{rank}_{\bc}(E,-J)+2c_1 (E,-J)\\
&=&\chi(\Sigma)\mathrm{rank}_{\bc}(E,J)-2c_1 (E,J),
\end{eqnarray*}
and finally, $H^1=\mathcal{H}=\mathcal{H}^+\oplus \mathcal{H}^-$ and
$$\kappa(\alpha,\alpha)=\int_\Sigma \Omega(\alpha\wedge
\star_{\Omega}\alpha)=||\alpha||^2$$ for $\alpha\in \mathcal{H}^+$, hence
positive definite, and similarly negative definite on $\mathcal{H}^-$.
Then
\begin{eqnarray*}
\mathrm{Signature}(\kappa)&=&\mathrm{dim}\mathcal{H}^+ -\mathrm{dim}\mathcal{H}^- \\
&=&\mathrm{Index}(d^-)-\mathrm{Index}(d^+)\\
&=&4c_1 (E,J).
\end{eqnarray*}
\end{pf}

\subsection{Link with Toledo invariants}

\begin{lemma}
\label{toledochern}
Let $\Gamma=\pi_{1}(\Sigma)$ be a surface group. For a symplectic linear representation $\rho:\Gamma\to Sp(V,\Omega)$,
\begin{itemize}
  \item the Toledo invariant equals the first Chern class of the flat symplectic vectorbundle $E_{\rho}$ associated to $\rho$.
  \item If $\rho(\Gamma)$ is contained in the subgroup that fixes a pair of complementary Lagrangian subspaces, then the Toledo invariant vanishes.
\end{itemize}\end{lemma}

\begin{pf}
Let $X(V,\Omega)$ denote the space of complex structures on $V$
tamed by $\Omega$ and compatible with $\Omega$. $Sp(V,\Omega)$ acts
by conjugation on $X(V,\Omega)$ transitively. Hence the metric on
$T_JX(V,\Omega)$ is determined by the killing form of
$\mathfrak{sp}(V,\Omega)$, which is proportional to $tr A^2$.
Tautologically, the trivial bundle $X(V,\Omega)\times V$ comes
equipped with an $Sp(V,\Omega)$-invariant complex structure $J$. The
trivial connection $D$ on it is not compatible with the complex
structure, but $\nabla=D-DJ$ is. This is an $Sp(V,\Omega)$-invariant
connection. Let $F=\frac{1}{2\pi}\mathrm{trace}(R^{\nabla})$ denote
the first Chern form obtained from its curvature. This is an
$Sp(V,\Omega)$-invariant 2-form on $X(V,\Omega)$, in fact, the
K\"ahler form of the symmetric metric on the Hermitian symmetric
space $X(V,\Omega)$ normalized so that the minimal sectional
curvature equals $-1$. Therefore, if $f:\tilde{\Sigma}\ra
X(V,\Omega)$ is a smooth $\rho$-equivariant map, the induced 2-form
$f^{*}F$ on $\Sigma$ is a representative of the first Chern class of
the flat symplectic vectorbundle $E_{\rho}$. Integrating it on
$\Sigma$ yields
\begin{eqnarray*}
T_{\rho}=c_{1}(E_{\rho}).
\end{eqnarray*}

Assume that $\rho$ fixes two complementary Lagrangian subspaces. In
other words, $V=R\oplus R^{*}$ is isomorphic to the sum of a real
vector space $R$ and its dual, with the tautological symplectic
structure $\Omega((q,p),(q',p'))=\langle p|q'\rangle-\langle
p'|q\rangle$. $\rho(\Gamma)$ falls into the subgroup $H\subset
Sp(V,\Omega)$ that preserves the splitting. This subgroup is
isomorphic to the full linear group $Gl(R)$. Fix a Euclidean
structure on $R$, whence an isomorphism $e:R\to R^{*}$. On $V$,
consider the complex structure $J=\begin{pmatrix}
0 &  e^{-1} \\
-e  & 0
\end{pmatrix}$. The tangent space of $X(V,\Omega)$ at $J$ consists in symmetric endomorphisms $A$ which anticommute with $J$.
For if $A_t$ is a curve in $X(V,\Omega)$ starting from $J$, the
differentiation gives $\dot A J+ J\dot A=0$. The complex structure
on this space is $\mathcal{J}A=JA$. The Euclidean structure on this
space is proportional to $\mathrm{Trace}(A^2)$. The tangent space
$W$ to the orbit $HJ$ is the subspace of symmetric endomorphisms $A$
which preserve the splitting $V=R\oplus R^{*}$. The elements $B=JA$
of $\mathcal{J}(W)$ exchange the summands. Therefore the products
$AB$ exchange the summands, $\mathrm{Trace}(AB)=0$, i.e. $W$ is
orthogonal to $\mathcal{J}(W)$. One concludes that the K\"ahler form
$F$ of $X(V,\Omega)$ vanishes on the orbit $HJ$. If
$f:\tilde{\Sigma}\ra X(V,\Omega)$ is a smooth $\rho$-equivariant
map, the induced 2-form $f^{*}F$ vanishes, thus $T_{\rho}=0$.
\end{pf}

\section{Appendix 2 : Structure of centralizers}

\begin{lemma}
\label{amencentr}
Let $G$ be a semisimple real algebraic group. Let $H\subset G$ be a semisimple subgroup with centralizer $\mathfrak{z}$. Then there exists an Iwasawa decomposition $G=KAN$ of $G$ such that $\mathfrak{z}$ is a direct sum of Lie algebras $\mathfrak{s}'\oplus\mathfrak{k}'\oplus\mathfrak{a}'$ where $\mathfrak{s}'$ is semisimple of noncompact type, $\mathfrak{k}'\subset \mathfrak{k}$ centralizes $\mathfrak{a}'$ and $\mathfrak{a}'\subset \mathfrak{a}$.
\end{lemma}

\begin{pf}
Since $H$ is semisimple, at least one of its orbits on $X=G/K$ is totally geodesic. All totally geodesic orbits of $H$ are pairwise equidistant. The union of these orbits is a totally geodesic subspace $Y$ of $X=G/K$, isometric to a Riemannian product $Y=W\times V\times F$, where $H$ acts on $W$ as an open subgroup of the full isometry group of $W$ and trivially on $V\times F$, $V$ is a symmetric space of noncompact type and $F$ is flat. Note that the isometry group $S'$ of $V$ canonically embeds in $G$, since it is generated by geodesic symmetries. Since it commutes with $H$, its Lie algebra $\mathfrak{s}'$ has to be contained in $\mathfrak{z}$. $\exp\mathfrak{z}$ leaves $Y$ invariant, thus it maps to $Isom(Y)$ with a compact kernel $K'$. Since $K'=Fix(Y)$ fixes $V$, it commutes with $S'$. Up to finite index, $Isom(Y)$ is a product $H\times A'$ where $A'$ is the group of translations along $F$. Since $A'$ commutes with $H$, its Lie algebra $\mathfrak{a}'\subset\mathfrak{z}$. Furthermore, $A'$ commutes with $K'$ and $S$. Fix points $w\in W$, $v\in V$ and $f\in F$, view $(w,v,f)$ as a point in $G$ and $F$ as a flat in $G$ passing through this point. Pick a maximal flat containing it. This defines an Iwasawa decomposition of $G$ such that $K'\subset K$ and $A'\subset A$. We conclude that $\mathfrak{z}=\mathfrak{s}'\oplus\mathfrak{k}'\oplus\mathfrak{a}'$, where $\mathfrak{k}'$ denotes the Lie algebra of $K'$.
\end{pf}

\begin{co}
\label{centralizerisreductive}
Let $G$ be a semisimple real algebraic group. Then centralizers of reductive subgroups of $G$ are reductive.
\end{co}

\begin{pf}
Let $H\subset G$ be reductive with centralizer $Z$. Write $H=SR$ where $S$ is semisimple and $R$ is the radical of $H$, a torus. Let $Z'$ be the centralizer of $S$. Then $R\subset Z'$ is a subtorus of $Z'$, $Z\subset Z'$ and $Z$ is the centralizer of $R$ in $Z'$. According to Corollary 2, section 13.17 of \cite{Borel_LAG}, in connected, affine algebraic groups, centralizers of subtori are reductive. Thus $Z$ is reductive.
\end{pf}

2000 {\sl{Mathematics Subject Classification.}}51M10, 57S25.

{\sl{Key words and phrases.}} Reductive surface group, flexible
representation, quadratic model, tube type Hermitian space \vskip1cm

\noindent     Inkang Kim\\
     School of Mathematics\\
     KIAS, Heogiro 87, Dongdaemen-gu\\
     Seoul, 130-722, Korea\\
     \texttt{inkang\char`\@ kias.re.kr}

\smallskip

     \noindent  Pierre Pansu\\ Laboratoire de Math{\'e}matiques
     d'Orsay\\
     UMR 8628 du CNRS\\
 Universit{\'e} Paris-Sud\\
 91405 Orsay C\'edex, France\\
  \texttt{pierre.pansu\char`\@ math.u-psud.fr}


\begin{thebibliography}{99}
{\small


\bibitem{Apa}Boris N. Apanasov, Bending and stamping deformations of hyperbolic manifolds.  Ann. Global Anal. Geom. {\bf 8} (1990), 3--12.

\bibitem{BGG} Steven B. Bradlow, Oscar Garc\'{\i}a-Prada, Peter B. Gothen, Surface group representations and ${\rm U}(p,q)$-Higgs bundles.
J. Differential Geom. \textbf{64}, (2003), 111--170.

\bibitem{Borel_coh}Armand Borel, Sur la cohomologie des espaces fibr\'es principaux et des espaces homog\`enes des groupes de Lie compacts. Ann. Math. {\bf 57},  (1953), 115--207.

\bibitem{Borel_LAG}Armand Borel, Linear algebraic groups. Second edition. Graduate Texts in Mathematics, 126. Springer-Verlag, New York, 1991.

\bibitem{BIW}Marc Burger, Alessandra Iozzi, Anna Wienhard, Surface
group representations with maximal Toledo invariant.
arXiv:math/0605656,  Annals of Math. \textbf{172}, (2010), no. 1,
517-566.

\bibitem{Eberlein} Patrick Eberlein; Geometry of nonpositively curved manifolds. Chicago Lectures in Mathematics. The University of Chicago Press, Chicago (1996).

\bibitem{G0} William Goldman, Discontinuous groups and the Euler class. Thesis, University of California at Berkeley (1980).

\bibitem{G1} William Goldman, Representations of fundamental groups of surfaces. Geometry and topology (J. Alexander and J. Harer, Eds), Lect. Notes Math. {\bf 1167}, Springer, 1985, 95--117.

\bibitem{GM} William Goldman and John Millson, Local rigidity of discrete groups acting on complex hyperbolic space. Invent. Math. {\bf 88}, (1987), 495--520.

\bibitem{GM3} William Goldman and John Millson, The deformation theory of representations of fundamental groups of compact K\" ahler manifolds. Publ. Math. I.H.\'E.S. {\bf 67}, (1988), 43--96.

\bibitem{HL} Luis Hern\`andez Lamoneda, Maximal representations of surface groups in bounded symmetric domains.
Trans. Amer. Math. Soc.  \textbf{324}, (1991), 405--420.

\bibitem{Hit1} Nigel Hitchin, The self-duality equations on a Riemann surface.  Proc. London Math. Soc. \textbf{55}, (1987), 59--126.

\bibitem{Hit2} Nigel Hitchin, Lie groups and Teichmüller space.  Topology  \textbf{31}, (1992), 449--473.

\bibitem{KPu} Inkang Kim and Pierre Pansu, Local rigidity
 in quaternionic hyperbolic space. JEMS \textbf{11}, (2009), 1141--1164.

\bibitem{classical} Inkang Kim and Pierre Pansu, Flexibility of surface groups in classical simple Lie groups. arXiv:math/

\bibitem{Manivel}Laurent Manivel, Fonctions sym\'etriques, polyn\^omes de Schubert et lieux de d\'eg\'enerescence. Cours Sp\'ecialis\'es, \textbf{3}. Soci\'et\'e Math\'ematique de France, Paris, (1998).

\bibitem{MDS}Dusa McDuff and Dietmar A. Salamon,
$J$-Holomorphic Curves and Symplectic Topology, Amer. Math. Soc. Colloquium Publications, Vol. \textbf{52}, (2004).

\bibitem{Meyer}Werner Meyer, Die Signatur von lokalen Koeffizientensystemen und Faserb\"undeln, Bonn. Math. Schr. \textbf{53} (1972).

\bibitem{Milnor}John Milnor, On the existence of a connection with curvature zero.  Comment. Math. Helv. \textbf{32}, (1958), 215--223.

\bibitem{NR}Albert Nijenhuis and Roger W. Richardson, Deformations of
homomorphisms of Lie groups and Lie algebras, Bull. Amer. Math. Soc.
\textbf{73}, (1967), 175--179.

\bibitem{Py}Pierre Py,
Indice de Maslov et th\'eor\`eme de Novikov-Wall.
Bol. Soc. Mat. Mexicana \textbf{11}, (2005), 303--331.

\bibitem{S} Carlos Simpson, Higgs bundles and local systems. Publ. Math. I.H.\'E.S. {\bf 75}, (1992), 5--95.

\bibitem{Th} William Thurston, The geometry and topology of 3-manifolds.
Lecture notes, Princeton, (1983).

\bibitem{To} Domingo Toledo, Representations of surface groups in complex
hyperbolic space, J. Diff. Geom. {\bf 29}, (1989), 125--133.}

\bibitem{Tu}Vladimir Turaev, A cocycle of the symplectic first Chern class and Maslov indices. Funktsional. Anal. i Prilozhen. \textbf{18}, (1984), 43--48.

    \end{thebibliography}
     \end{document}